\documentclass[final,11pt,english]{scrartcl}
\usepackage{amsmath,amssymb,amsfonts,amsthm,mathtools}
\usepackage[hidelinks,pdfusetitle]{hyperref}
\usepackage{babel}
\usepackage[utf8]{inputenc}
\usepackage[T1]{fontenc}
\usepackage{microtype,lmodern}
\usepackage{fullpage}

\usepackage{ifdraft}
\usepackage[inline]{showlabels}

\mathtoolsset{mathic}

\usepackage[
		backend=biber,
		url=false,
		style=alphabetic,
	]{biblatex}
\usepackage{csquotes}
\usepackage[noabbrev,capitalize]{cleveref}

\usepackage[obeyFinal]{todonotes}
\usepackage{easyReview}

\usepackage{multifootnote}
\usepackage{nicematrix}
\usepackage{cancel}
\usepackage{verbatim}

\usepackage{tikz}
\usetikzlibrary{babel,cd}

\newcommand{\email}[1]{\href{mailto:#1@uni-hamburg.de}{#1@uni-hamburg.de}}

\author{Vicente Cortés\texorpdfstring{\multifootnotemark[fn:vc]}{} \and %
Shubham Dwivedi\texorpdfstring{\multifootnotemark[fn:sd]}{} \and %
Lorenz Jetter\texorpdfstring{\multifootnotemark[fn:lj]}{}}
\title{Generalized Ricci Solitons on Heisenberg and Oscillator Algebras}
\date{November 1, 2026}

\hypersetup{
	pdfkeywords = {Generalized Ricci Solitons, Generalized Algebraic Solitons, Nilmanifolds, Heisenberg Algebras, Oscillator Algebras}
}

\usepackage[inline]{enumitem}
\setlist[enumerate,1]{label=(\roman*)}

\usepackage{diffcoeff}%

\usepackage{tensor}

\usepackage{lmacros} %

\usepackage{blindtext}

\newcommand{\heis}{\mathfrak{h}}
\DeclareMathOperator{\Aut}{Aut}
\DeclareMathOperator{\Der}{Der}
\DeclareMathOperator{\Sp}{Sp}
\DeclareMathOperator{\GL}{GL}

\DeclareMathOperator{\OO}{O}
\DeclareMathOperator{\SO}{SO}
\DeclareMathOperator{\UU}{U}
\DeclareMathOperator{\SU}{SU}

\newcommand{\GG}{\mathrm{G}}
\newcommand{\NN}{\mathrm{N}}

\newcommand{\frakG}{\mathfrak{G}}
\newcommand{\bmu}{{\boldsymbol{\mu}}}

\newcommand{\frakn}{\mathfrak{n}}
\newcommand{\frakN}{\mathfrak{N}}

\DeclareMathOperator{\Diag}{Diag}

\providecommand{\longhookrightarrow}{\lhook\joinrel\longrightarrow}

\begin{document}

\maketitle

\multifootnotetext[fn:vc,fn:sd,fn:lj]{%
		\textsc{Universität Hamburg}, Fachbereich Mathematik, %
		Bundesstraße~55, 20146 Hamburg, Germany \\
		E-mail: \email{vicente.cortes}, \email{shubham.dwivedi}, \email{lorenz.jetter}}

\begin{abstract}
We prove that on a connected, simply-connected nilpotent Lie group, the notions of left-invariant generalized Ricci soliton, semi-algebraic and algebraic generalized Ricci soliton coincide.
Moreover, every such soliton has harmonic torsion.
We classify all algebraic solitons with positive definite metric on the five-dimensional Heisenberg algebra, the four-dimensional oscillator algebra and the four-dimensional split-oscillator algebra up to automorphisms and rescalings.
We find that there is only a finite number of algebraic solitons, all of which have harmonic torsion.
\bigskip

\smallskip

\end{abstract}

\tableofcontents

\makeatletter
\let\@oldthetheorem\thetheorem
\renewcommand{\thetheorem}{\Alph{theorem}}
\makeatother

\section{Introduction}

A generalized algebraic soliton on a finite-dimensional, real Lie algebra~$\frakg$ is a triple~$(g, H, \lambda)$ with $g \in \Sym_+^2 \frakg^*$ a positive definite, symmetric bilinear form, $\lambda \in \R$ and $H \in \Lambda^3_{\mathrm{cl}} \frakg^*$ a closed three-form on $\frakg$, solving
\begin{equation} \label{eq:nil_soliton_intro}
	\left\{
	\begin{aligned}
		\Ric^B_{\mu, H} &= \lambda \Id + D, \\
		-\D \D^*_g H &= \lambda H + \rho\pa[\Big]{\Ric^B_{\mu, H}} H,
	\end{aligned}
	\right.
\end{equation}
where $D \in \Der(\frakg)$ is a derivation of the Lie bracket~$\mu$ of $\frakg$, $\D^*_g$ is the co-differential with respect to $g$,
\begin{inline}
	\rho\colon \mathfrak{gl}(\frakg) \to \End\pa*{\Lambda^3\frakg^*}
\end{inline}
is the canonical representation of $\mathfrak{gl}(\frakg)$ on three-forms and $\Ric^B_{\mu, H} = g^{-1} \pa{\Ric^g - \frac{1}{4} H_g^2} \in \mathfrak{gl}(\frakg)$ is an endomorphism with $H_g^2(X, Y) \ceq g(X\lrcorner H, Y\lrcorner H)$ and $\Ric^g$ the Ricci-tensor induced by $g$.
Introduced in \cite{Fusi2024}, this notion generalizes that of an algebraic soliton of the bracket flow \cite{Lauret2001}, the latter of which corresponds to the special case of $H = 0$ in \cref{eq:nil_soliton_intro}.
Generalized algebraic solitons are equivalent to special solutions of the generalized bracket flow in the context of Hitchin's generalized geometry and provide important examples of generalized Ricci solitons.
More precisely, an algebraic soliton may be understood as flow of a Dorfman bracket on a left-invariant, exact, metric Courant algebroid over a Lie group~$\GG$ which evolves only by scaling.

A natural question is whether the algebraic soliton equations capture all left-invariant self-similar solutions of the generalized Ricci flow.
A generalized Ricci soliton is initially defined using diffeomorphisms, which need not preserve the underlying Lie group structure.
Thus, left-invariance of the metric and torsion does not by itself identify the soliton vector field with a derivation of the Lie algebra.
For ordinary Ricci solitons on simply-connected nilpotent Lie groups, the passage from self-similarity to algebraic soliton equations is fundamental to the theory of nilsolitons \cite{Lauret2001,Jablonski2014}.
Our first main result establishes the corresponding statement for the generalized Ricci flow.

\begin{theorem}[{\cref{thm:nilpotent-geometric-algebraic}}]
Let $\NN$ be a connected, simply-connected nilpotent Lie group, let $g$ be a left-invariant Riemannian metric and let $H$ be a left-invariant closed three-form.
Then $(g, H, \lambda)$ is a generalized Ricci soliton if and only if there exists a $g$-symmetric derivation $D \in \Der(\mathfrak n)$, where $\frakn \ceq \operatorname{Lie}(\NN)$, such that
\begin{align*}
	g^{-1}\left({\Ric}^{g}-\frac14H_g^2\right) =\lambda \Id+D, \qquad \D_g^*H=0, \qquad \rho(D)H=2\lambda H.
\end{align*}
In particular, the notions of generalized Ricci soliton, semi-algebraic and algebraic generalized Ricci soliton coincide in this setting.
\end{theorem}
The theorem applies to the original left-invariant pair~$(g, H)$, without assuming that its self-similarity is generated by automorphisms or that its torsion is co-closed. It therefore reduces the classification of left-invariant geometric generalized Ricci solitons on simply-connected nilpotent Lie groups to a finite-dimensional algebraic problem.
\medbreak

In view of the Levi-Malcev decomposition, generalized solitons on nilpotent or solvable Lie groups form a particularly important class of left-invariant generalized solitons on Lie groups.
In \cite{Fusi2024}, all generalized algebraic solitons over nilpotent Lie algebras were classified in dimension three and four, in particular including the three-dimensional Heisenberg algebra (as the unique, three-dimensional, non-abelian, nilpotent Lie algebra).
We extend this result in \cref{sec:heisenberg_algebras} by classifying all generalized algebraic solitons on the five-dimensional Heisenberg algebra~$\heis_5$ up to automorphism and rescaling (in the sense of \cref{lem:rescaling_soliton}).

\begin{theorem}[{\cref{thm:grs_heis5}}] \label{thm:intro_grs_heis5}
Up to automorphism and rescaling, all nilsolitons~$(g, H, \lambda)$ on the five-dimensional Heisenberg algebra in the ordered basis~$\pa{p_1, p_2, q_1, q_2, z}$ are given by $[g_{ij}] = [\delta_{ij}]$, $\lambda = -2$ and either
\begin{enumerate}
\item $H \ceq 0$ (classical),
\item $H \ceq \pa*{p^1 \wedge q^1 - p^2 \wedge q^2} \wedge z^*$,
\item $H \ceq \pa*{p^1 \wedge p^2 - q^1 \wedge q^2} \wedge z^*$ or
\item $H \ceq \frac{1}{3} \sqrt{3} \cdot \pa*{p^1 \wedge p^2
			+ p^1 \wedge q^1
			- p^2 \wedge q^2
			- q^1 \wedge q^2} \wedge z^*$,
\end{enumerate}
where the Lie bracket is given by $[p_1, q_1] = [p_2, q_2] = z$ and zero otherwise.
None of the above solitons are isomorphic.
\end{theorem}

Note, that due to \cref{thm:nilpotent-geometric-algebraic}, we in particular classified all left-invariant generalized Ricci solitons on the Heisenberg group~$H_5$.
The proof proceeds by first showing in \cref{lem:simplification_metric} that all positive definite, symmetric bilinear forms~$g \in \Sym_+^2 \heis_n^*$ on an $n$-dimensional Heisenberg algebra~$\heis_n$ (with $n = 2k + 1$) can be transformed to
\begin{inline}
	[g\indices{_{ij}}]
	=
		\Diag(D, D, 1)
\end{inline}
with $D \in \R^{k \times k}$ a positive definite diagonal matrix.
Then, we show that necessarily, nilsolitons~$(g, H, \lambda)$ on $\heis_5$ have harmonic torsion.
\begin{lemma}[{\cref{lem:heis5_solitons_harmonic}}] %
All nilsolitons on the five-dimensional Heisenberg algebra have harmonic torsion.
\end{lemma}
We do so by direct solution of the soliton equations after the above-mentioned simplification of the metric, thereby obtaining $\Laplace^g H = 0$ as a necessary constraint.
In the more recent paper \cite{Chen2026}, it was shown that all generalized nilsolitons in general have harmonic torsion and in \cite{Fusi2026}, this result is further extended to generalized algebraic solitons on simply-connected, unimodular Lie groups with positive semi-definite Killing form using more involved arguments our alternative, direct proof for $\heis_5$ does not rely on.
The last part of the proof of \cref{thm:intro_grs_heis5} consists in studying the action of the stabilizer~$\Aut(\heis_5)_g = U(2) \rtimes \Z_2$ of the metric on all harmonic three-forms on $\heis_5$, which reduces the nilsoliton equations to a pair of simple equations in two variables.
\smallbreak

In \cref{sec:oscillator_algebras}, we classify all generalized algebraic solitons on the four-dimensional oscillator algebra~$\mathfrak{osc}_4^{-1}$ and the four-dimensional split-oscillator algebra~$\mathfrak{osc}_4^1$, yielding the following results:
\begin{theorem}[{\cref{thm:grs_osc}}] \label{thm:intro_grs_osc}
There do not exist any algebraic solitons (with positive definite background metric~$g$) on the four-dimensional oscillator algebra.
\end{theorem}

\begin{theorem}[{\cref{thm:grs_split_osc}}] \label{thm:intro_grs_split_osc}%
Up to automorphism and rescaling of the metric, all algebraic solitons~$(g, H, \lambda)$ on the four-dimensional hyperbolic oscillator algebra in the ordered basis~$\pa{d, p, q, z}$ are given by $\lambda = -2$
and either
\begin{enumerate}
\item
\begin{inline}
	[g_{ij}]
		=
			\Diag\pa*{1,\, \frac{1}{2} \sqrt{3},\, \frac{1}{2} \sqrt{3},\, 1},
\end{inline}
$H = 0$ (classical),
\item $[g_{ij}] = [\delta_{ij}]$, $H = +p^* \wedge q^* \wedge z^*$ or
\item $[g_{ij}] = [\delta_{ij}]$, $H = -p^* \wedge q^* \wedge z^*$,
\end{enumerate}
where the Lie bracket is given by
\begin{inline}
	[d, p] = q,\, [d, q] = p,\, [p, q] = z
\end{inline}
and zero otherwise.
None of these three solitons are isomorphic.
\end{theorem}
Independently of our paper, the results of \cref{thm:intro_grs_osc,thm:intro_grs_split_osc} were also recently discovered in \cite{Mucciante2026}, where the author extended the classification of generalized algebraic solitons on nilpotent Lie algebras in \cite{Fusi2024} to all solvable algebras in dimension three and four, also establishing that on solvable Lie algebras, all generalized algebraic solitons with non-vanishing torsion are expanding.
The proof in \cite{Mucciante2026} uses the classification of left-invariant metrics on four-dimensional, simply-connected, unimodular Lie groups \cite{Thuong2017}, while our proof for both theorems is self-contained and works by first simplifying a generic metric~$g \in \Sym_+^2 (\mathfrak{osc}_4^{\tau})^*$ with $\tau \in \set{\pm 1}$ analogous to the proof in the Heisenberg case (and treating both algebras simultaneously).
We then determine all solutions of the solvsoliton equations through several big case distinctions.

\makeatletter
\let\thetheorem\@oldthetheorem
\makeatother

\paragraph*{Conventions}

If $\set{e_i}$ is an orthonormal base of an $n$-dimensional Lie algebra~$\frakg$ with respect to some positive definite, symmetric bilinear form~$g \in \Sym^2_+ \frakg^*$, the inner product on $(0, k)$-tensors is defined by requiring all $\iset{e^{j_1} \tprod \dotsm \tprod e^{j_k}}{\range{j_1}{j_k} \in \setrange{1}{n}}$ to be orthonormal, where $\set{e^j}$ is the dual-base of $\set{e_i}$.
The Laplacian induced by $g$ is set to be
\begin{inline}
	\Laplace_g \ceq -\pa{\D \D^*_g + \D^*_g \D}
\end{inline}
with the Chevalley-Eilenberg differential~${\D}$ and the co-differential~$\D^*_g = (-1)^k \hodge^{-1} \D\, *$ on $k$-forms, where the Hodge star satisfies
\begin{inline}
	\alpha \wedge \hodge\beta = \frac{1}{k!} g(\alpha, \beta) \vol^g
\end{inline}
for all $\alpha, \beta \in \Lambda^k\frakg^*$.

\paragraph*{Acknowledgments}
This work is funded by the Deutsche Forschungsgemeinschaft (DFG, German Research Foundation) under Germany's Excellence Strategy -- EXC~2121 ``Quantum Universe'' -- 390833306.
VC is also funded by SFB-Geschäftszeichen~1624 -- Projektnummer 506632645.
The authors would like to thank James Stanfield and Lorenzo Mucciante for their insightful comments and ongoing discussions.

\section{Preliminaries}

In this section, we briefly repeat well-known results about the generalized Ricci flow.
We then specialize on left-invariant solutions on Lie groups (or more generally $\GG/\Gamma$ with $\Gamma \subset \GG$ a lattice in $\GG$ acting from the right) and introduce the generalized bracket flow.
In then end, we state the definition of generalized algebraic solitons and some important properties thereof.
\smallbreak

In the following, $M$ will always denote a smooth manifold.
Given an arbitrary Riemannian metric~$g$ on $M$ and a closed three-form~$H \in \Omega_{\mathrm{cl}}^3(M)$, there exists a unique metric connection~$\nabla^{g, H}$ with torsion equal to $H$, known as the associated \emph{Bismut connection}.
Denoting the Levi-Civita connection of $g$ by $\nabla^g$, it is determined by
\begin{equation*}
	g\pa*{\nabla^{g, H}_X Y, Z}
	= g\pa*{\nabla^g_X Y, Z} + \frac{1}{2} H(X, Y, Z).
\end{equation*}
One verifies that the Ricci curvature of the Bismut connection is
\begin{equation} \label{eq:ricci_bismut}
	\Ric^{g, H} = \Ric^g - \frac{1}{4} H^2 - \frac{1}{2} \D^*_g H,
\end{equation}
where $H^2$ is the symmetric $(0, 2)$-tensor given by $H^2(X, X) \ceq \abs{X \lrcorner H}_g^2$.
In local coordinates,
\begin{equation*}
	{(H)^2}\indices{_{ij}}
	= \tensor{g}{^{ac}} \tensor{g}{^{bd}} \tensor{H}{_{abi}} \tensor{H}{_{cdj}}.
\end{equation*}
All computational details hereof may be found in \cite[Section~3.3]{Garcia-Fernandez2021}.
Following the notation in \cite{Fusi2024}, we denote by $\Ric^B_{g, H}$ the \emph{symmetric part} of the Bismut Ricci curvature~$\Ric^{g, H}$, that is,
\begin{equation*}
	\Ric^B_{g, H} = \Ric^g - \frac{1}{4} H^2.
\end{equation*}
Trivially, $\Ric^B_{g, H}$ and $\Ric^{g, H}$ agree if and only if $H$ is harmonic.

Given an initial metric~$g_0$ and closed~$H_0$, we recall \cite[Lemma~4.6]{Garcia-Fernandez2021} the following weakly parabolic flow:
\begin{definition}[Generalized Ricci Flow]
A smooth one-parameter family~$(g(t), b(t))$ consisting of Riemannian metrics~$g(t)$ and two-forms~$b(t)$ solves the \emph{generalized Ricci flow} equations with initial value~$(g_0, H_0)$, if
\begin{equation} \label{eq:grf_bismut}
	\begin{cases}
		\partial_t\pa*{g(t) - b(t)} = -2 \Ric^{g(t), H_0 + \D b(t)}, \\
		g(0) - b(0) = g_0.
	\end{cases}
\end{equation}
\end{definition}

This flow is a generalization of the Ricci flow, which corresponds to solutions of \cref{eq:grf_bismut} with vanishing initial torsion, $H_0 = 0$.
Importantly, \cref{eq:grf_bismut} is equivalent to the flow
\begin{equation*}
	\calG^{-1} \diffp*{\calG}{t} = -2 \mathcal{R}c(\calG),\qquad
	\calG(0) = \calG_0
\end{equation*}
of so-called ``generalized metrics'' in the setting of Hitchin's generalized geometry.
In the following, we will elaborate on this correspondence and briefly introduce the most important notions of generalized geometry needed for understanding the results in this paper:
\begin{definition}[Courant Algebroid]
A Courant algebroid~$(E, \pi, [\cdot, \cdot], \ipr{\cdot}{\cdot})$ over a manifold~$M$ consists of
\begin{itemize}
\item A vector bundle $E \to M$,
\item A smooth vector bundle homomorphism~$\pi\colon E \to TM$, called \emph{anchor map},
\item An $\R$-bilinear map~$[\cdot, \cdot]\colon \Gamma(E) \times \Gamma(E) \to \Gamma(E)$, called \emph{Dorfman bracket} of $E$ and
\item A symmetric, non-degenerate bilinear form~$\ipr{\cdot}{\cdot} \in \Gamma(\Sym^2 E^*)$, called \emph{neutral inner product} of $E$,
\end{itemize}
satisfying the following three axioms for all $a, b, c \in \Gamma(E)$:
\begin{enumerate}[series=caconditions]
\item The Dorfman bracket is a derivation of itself, that is,
\begin{inline}
	\bra*{a, [b, c]}
	= \bra*{[a, b], c} + \bra*{b, [a, c]}.
\end{inline}
\item The bracket is compatible with the inner product, in the sense that
\begin{equation*}
	\pi(a)\ipr{b}{c}
	= \ipr{[a, b]}{c} + \ipr{b}{[a, c]}.
\end{equation*}
\item \label{it:ca_axiom_symm_part}
The symmetric part of the bracket satisfies
\begin{inline}
	\ipr{[u, v] + [v, u]}{w} = \pi(w)\ipr{u}{v}.
\end{inline}
\end{enumerate}

Two Courant algebroids~$(E_1, \pi_1, [\cdot, \cdot]_1, \ipr{\cdot}{\cdot}_1)$ and $(E_2, \pi_2, [\cdot, \cdot]_2, \ipr{\cdot}{\cdot}_2)$ over $M$ are isomorphic if there exists a vector bundle isomorphism~$F\colon E_1 \to E_2$ which preserves the extra structures of the algebroids, i.e.\
\begin{inline}
	\pi_2 \circ F = \pi_1,\, F\bra*{\cdot, \cdot}_1 = \bra*{F \cdot, F \cdot}_2\ \text{and}\ F^*\ipr{\cdot}{\cdot}_2 = \ipr{\cdot}{\cdot}_1.
\end{inline}
\end{definition}

It is easy to prove that every Courant algebroid additionally satisfies the following properties for sections $a, b \in \Gamma(E)$ and $f \in \calC^\infty(M)$:
\begin{enumerate}[resume=caconditions]
\item The anchor is compatible with the Dorfman bracket, $\pi [a, b] = [\pi a, \pi b]$,
\item The Dorfman bracket satisfies the \emph{anchored Leibniz rule},
\begin{inline}
	[a, fb] = L_{\pi(a)}(f)\, b + f[a, b].
\end{inline}
\end{enumerate}

The map~$\pi^*\colon T^*M \to E$ stands for the dual~$T^*M \to E^*$ of the anchor~$\pi\colon E \to TM$, where we metrically identify $E$ with $E^*$.
It is known that on every Courant algebroid,
\begin{inline}
	\ipr{\pi^* \cdot}{\pi^* \cdot} = 0
\end{inline}
or, stated differently, that
\begin{equation} \label{eq:eca_short_sequence}
\begin{tikzcd}
	0 \ar[r] & T^*M \ar[r,"\pi^*"] & E \ar[r, "\pi"] & TM \ar[r] \ar[l, bend left, "\sigma"] & 0
\end{tikzcd}
\end{equation}
constitutes a well-defined short sequence.
One verifies that furthermore,
\begin{inline}
	\bra*{\pi^* \cdot, \pi^* \cdot} = 0.
\end{inline}

In this paper, we only consider \emph{exact} Courant algebroids~(ECAs), that is,
Courant algebroids for which sequence~\labelcref{eq:eca_short_sequence} is exact.

As we shortly review now, it is a fundamental fact of the theory of Courant algebroids due to \v{S}evera, that every exact Courant algebroid is isomorphic to a twisted generalized tangent bundle:
\begin{definition}[Twisted Generalized Tangent Bundle]
For a smooth manifold~$M$ and three-form~$H \in \Omega^3(M)$, consider the vector bundle
\begin{inline}
	\mathbb{T}M \ceq TM \oplus T^*M
\end{inline}
together with
\begin{gather*}
	\pi\colon \mathbb{T}M \longto TM,\quad X + \xi \longmapsto X, \\
	\ipr{\cdot}{\cdot} \in \Gamma(\Sym^2 \mathbb{T}M^*),\quad \ipr{X + \xi}{Y + \eta} \ceq \frac{1}{2}\pa*{\xi(Y) + \eta(X)}
	\shortintertext{and}
	[\cdot, \cdot]_H\colon \Gamma(\mathbb{T}M) \times \Gamma(\mathbb{T}M) \longto \Gamma(\mathbb{T}M),\\
			[X + \xi, Y + \eta]_H \ceq L_X Y + L_X \eta - \iota_Y \D\xi + \iota_Y \iota_X H
\end{gather*}
for all $X, Y \in \VF(M)$ and $\xi, \eta \in \Omega^1(M)$.
Then the data~$\pa*{\mathbb{T}M, \pi, \bra*{\cdot, \cdot}_H, \ipr{\cdot}{\cdot}}$ defines an exact Courant algebroid if and only if $H$ is closed, which we call the \emph{$H$-twisted generalized tangent bundle} of $M$.
Similarly, $\bra*{\cdot, \cdot}_H$ is referred to as the \emph{($H$-twisted) Dorfman bracket}.
In the special case of $H \equiv 0$, the generalized tangent bundle is said to be \emph{untwisted}.
\end{definition}

Every exact Courant algebroid admits an isotropic splitting~$\sigma\colon TM \to E$ of the short exact sequence in \labelcref{eq:eca_short_sequence}. %
As a consequence, $E = \Image \pi^* \oplus \Image \sigma \cong T^*M \oplus TM$ is a decomposition of $E$ into isotropic sub-bundles.
In fact, one can show that
\begin{inline}
	H_\sigma \ceq 2\ipr{[\sigma \cdot, \sigma \cdot]}{\sigma\cdot} \in \Omega^3(M)
\end{inline}
defines a closed three-form and $F_\sigma\colon TM \oplus T^*M \to E$, $F_\sigma(X + \xi) \ceq \sigma X + \frac{1}{2}\pi^* \xi$ an isomorphism between $\mathbb{T}M_{H_\sigma}$ and $E$.
Thus, we obtain the following correspondence:
\begin{theorem}[\v{S}evera]
Let $\sigma$ and $\tilde{\sigma}$ be two isotropic splittings of the exact Courant algebroid~$E$.
Then the induced three-forms~$H_\sigma$ and $H_{\tilde{\sigma}}$ lie in the same de~Rham cohomology class.
In particular, there is a one-to-one correspondence between $H^3_{\mathrm{dR}}(M)$ and the isomorphism classes of exact Courant algebroids over $M$.
\end{theorem}
\medbreak

Next, we briefly state the classification of generalized metrics on exact Courant algebroids:
\begin{definition}
Let $\OO(E, \ipr{\cdot}{\cdot})$ be all Courant algebroid automorphisms which are orthogonal with respect to the inner product~$\ipr{\cdot}{\cdot}$.
A \emph{generalized metric~$\calG$} on an ECA~$(E, \pi, \ipr{\cdot}{\cdot}, [\cdot, \cdot])$ is a vector bundle automorphism~$\calG \in \OO(E, \ipr{\cdot}{\cdot})$, such that $\ipr{\calG \cdot}{\cdot}$ is symmetric and positive definite.
An exact Courant algebroid equipped with a generalized metric is referred to as \emph{metric (exact) Courant algebroid}.
\end{definition}

Given a generalized metric~$\calG$ on $E$, there exists a unique isotropic splitting~$\sigma$, such that on $TM \oplus T^*M$, the metric reads
\begin{equation} \label{eq:generalized_metric_correspondence}
	F_\sigma^{-1} \calG F_\sigma
	=
			\begin{bmatrix}
				0 & g^{-1} \\
				g & 0
			\end{bmatrix}
\end{equation}
with $g$ a Riemannian metric on $M$ (see \cite[Lemma~2.29]{Garcia-Fernandez2021} for details).
This splitting is said to be \emph{adapted to $\calG$} or \emph{induced by $\calG$}.
Conversely, the expression in \cref{eq:generalized_metric_correspondence} defines a generalized metric on $E$, which we call~$\calG(g, \sigma)$.

If $\sigma, \tilde{\sigma}\colon TM \to E$ are two isotropic splittings of $E$, then
\begin{equation*}
	b \ceq \pa*{\pi^*}^{-1} \circ \pa*{\tilde{\sigma} - \sigma}
	= \ipr{\tilde{\sigma} \cdot}{\sigma \cdot} \in \Omega^2(M)
\end{equation*}
is a well-defined two-form.
On the level of the generalized tangent bundle, the diagram
\begin{equation*}
	\begin{tikzcd}
		E \ar[r, "F_\sigma"] \ar[dr, "F_{\tilde{\sigma}}"below]
				& \mathbb{T}M_{H} \ar[d, "e^b", "\cong"left] \\
		& \mathbb{T}M_{H + \D b}
	\end{tikzcd}\qquad
	\text{with}\quad e^b
	\ceq
			\begin{bmatrix}
				\Id & 0 \\
				b & \Id
			\end{bmatrix}
\end{equation*}
commutes.
Every metric on the generalized tangent bundle~$\mathbb{T}M_H$ is therefore given by
\begin{equation} \label{eq:gen_metric_TTM}
	\calG(g, b) \ceq e^b \begin{bmatrix}
		0 & g^{-1} \\
		g & 0
	\end{bmatrix} e^{-b}
\end{equation}
for a uniquely determined Riemannian metric~$g$ and two-form~$b \in \Omega^2(M)$.
\bigbreak

In order to define the generalized Ricci curvature, we briefly introduce generalized connections and divergence operators on exact Courant algebroids:

\begin{definition}[Generalized Connections]
A \emph{generalized connection} is a map
\begin{equation*}
	D\colon \Gamma(E) \longto \Gamma(E^* \tprod E),\qquad
	D_u v \ceq (Dv)(u)
\end{equation*}
which satisfies the anchored Leibniz rule
\begin{inline}
	D_u(fv) = L_{\pi(u)}(f)\, v + f D_u v
\end{inline}
and is compatible with the inner product, i.e.\
\begin{inline}
	L_{\pi(u)}\ipr{v}{w}
	= \ipr{D_u v}{w} + \ipr{v}{D_u w}
\end{inline}
for all $u, v, w \in \Gamma(E)$ and $f \in \calC^\infty(M)$.
We furthermore define the \emph{torsion} of a generalized connection as
\begin{equation*}
	T^D(u, v)
	\ceq D_u v - D_v u - [u, v] + (Du)^* v
\end{equation*}
for $u, v \in \Gamma(E)$.
\end{definition}

\begin{definition}[Generalized Divergence]
A \emph{generalized divergence operator} is a map
\begin{equation*}
	\Div\colon \Gamma(E) \longto \calC^\infty(M)
\end{equation*}
which satisfies the anchored Leibniz rule.
That is,
\begin{inline}
	\Div(f u) = L_{\pi u}(f) + f \Div u
\end{inline}
for all $u \in \Gamma(E)$ and $f \in \calC^\infty(M)$.
\end{definition}
We define the following canonical induced divergence operators:
\begin{enumerate}
\item For a generalized connection~$D$, we set
\begin{inline}
	\Div^D(u) \ceq \Trace Du
\end{inline}
and
\item For a generalized metric~$\calG$,
\begin{inline}
	\Div^\calG(u) \ceq \Div^g(\pi u),
\end{inline}
where $g$ is the Riemannian metric induced by $\calG$ and $\Div^g\colon \VF(M) \to \calC^\infty(M)$ the divergence on $M$.
\end{enumerate}

\begin{definition}[Generalized Riemann Curvature]
The \emph{generalized Riemann curvature} associated to a torsion free generalized connection is the tensor
\begin{inline}
	{\mathcal{R}m}^D \in \Gamma\pa*{\Lambda^2 E^* \tprod \End E}
\end{inline}
which is defined via
\begin{equation*}
	\ipr{{\mathcal{R}m}^D(v, w)b}{a}
	\ceq \frac{1}{2}\bra*{\ipr*{\pa[\big]{D_{v, w}^2 - D_{w, v}^2} b}{a} + \ipr*{\pa[\big]{D_{b, a}^2 - D_{a, b}^2} v}{w}
			- \ipr*{(Dv)^*w}{(Db)^* a}}.
\end{equation*}
By taking the trace, we get the \emph{full generalized Ricci curvature}
\begin{inline}
	\overline{\mathcal{R}c}^D(u, v)
	\ceq \Trace {\mathcal{R}m}^D(\cdot, u) v.
\end{inline}
\end{definition}

\begin{definition}
A generalized connection is said to be \emph{Levi-Civita} with respect to a generalized metric~$\calG$, if $\calG$ is parallel, that is,
\begin{inline}
	D(\calG) \ceq D \circ \calG - \calG \circ D = 0,
\end{inline}
and the torsion~$T^D$ vanishes.
\end{definition}

\begin{definition}[Generalized Ricci Curvature]
Let $E = E_+ \oplus E_-$ be the decomposition of $E$ into the eigenbundles of $\calG$ (to the eigenvalues~$\pm 1$).
For a generalized Levi-Civita connection~$D$, we consider the restrictions of the full generalized Ricci tensor to the eigenbundles
\begin{equation*}
	{\mathcal{R}c}^+ \ceq \overline{\mathcal{R}c}^D|_{E_- \times E_+}
\qquad\text{and}\qquad
	{\mathcal{R}c}^- \ceq \overline{\mathcal{R}c}^D|_{E_+ \times E_-}.
\end{equation*}
The \emph{generalized Ricci tensor} of $D$ is then defined as the endomorphism
\begin{inline}
	{\mathcal{R}c}^D \in \End(E)
\end{inline}
obtained by raising one index of ${\mathcal{R}c}^+ + {\mathcal{R}c}^-$ via $\calG$, that is,
\begin{inline}
	\ipr{\calG\, {\mathcal{R}c}^D \cdot}{\cdot} = {\mathcal{R}c}^+ + {\mathcal{R}c}^-.
\end{inline}
\end{definition}
Unlike in the non-generalized setting, for a fixed generalized metric~$\calG$, the Levi-Civita conditions do not uniquely determine the connection.
However, it is known that there always exists a generalized Levi-Civita connection with prescribed divergence~$\Div$ and that for two Levi-Civita connections~$D$ and $\tilde{D}$ with same divergence, the generalized Ricci curvature agrees, ${\mathcal{R}c}^D = {\mathcal{R}c}^{\tilde{D}}$.
It is therefore well-defined to denote this induced Ricci tensor by ${\mathcal{R}c}(\calG, \Div)$.

For the purposes of this paper, it is enough to consider generalized Levi-Civita connections with metric divergence~$\Div^\calG$.
For $\calG = \calG(g, 0)$ on ${\mathbb{T}M}_H$, the generalized Ricci curvature takes the form \cite[Prop.~3.30 + Def.~3.31]{Garcia-Fernandez2021}
\begin{equation} \label{eq:gen_ricci_gen_tangent_bundle}
	\mathcal{R}c(\calG)
	\ceq \mathcal{R}c(\calG, \Div^\calG)
	=
			\begin{bmatrix}
				g^{-1} \Ric_{g, H}^B & \frac{1}{2} g^{-1} \D^*_g H g^{-1} \\
				-\frac{1}{2} \D^*_g H & -\Ric_{g, H}^B g^{-1}
			\end{bmatrix}.
\end{equation}

\subsection{Generalized Ricci Solitons}

Let
\begin{equation*}
	\calM^E \ceq \iset{\calG \in \OO(E, \ipr{\cdot}{\cdot})}{\calG^2 = \Id_E\ \text{and}\ \ipr{\calG\cdot}{\cdot}\ \text{is pos.\ def.}}
\end{equation*}
denote the sheaf of all generalized metrics on an ECA~$(E, \pi, \ipr{\cdot}{\cdot} [\cdot, \cdot])$.
The tangent space of a generalized metric~$\calG \in \calM^E$ is given by \cite[Lemma~2.31]{Garcia-Fernandez2021}
\begin{equation*}
	T_\calG\calM^E
	= \iset*{\calG\Phi \in \End(E)}{\Psi \in \mathfrak{so}(E, \ipr{\cdot}{\cdot}),\, \Phi \calG + \calG \Phi = 0}.
\end{equation*}
For a smooth one-parameter family of generalized metrics~$(\calG(t))_{t \in I}$ with $\calG_0 = \calG$, this in particular implies
\begin{inline}
	\calG \dot{\calG} \in \mathfrak{so}(E, \ipr{\cdot}{\cdot})
\end{inline}
and under the isomorphism~$(E, \calG_0) \cong \pa*{(T \oplus T^*)_H, \calG(g_0, 0)}$, it is
\begin{equation*}
	\calG\dot{\calG}
	=
			\begin{bmatrix}
				g^{-1} \dot{g} & g^{-1} \dot{b} g^{-1} \\
				-\dot{b} & -\dot{g} g^{-1}
			\end{bmatrix}.
\end{equation*}

Together with expression~\labelcref{eq:gen_ricci_gen_tangent_bundle}, we observe that the generalized Ricci flow defined in \cref{eq:grf_bismut} and the flow
\begin{equation} \label{eq:grf}
	\calG^{-1} \diffp*{\calG}{t} = -2 \mathcal{R}c(\calG)
\end{equation}
of generalized metrics~$\calG(t)$ with initial value~$\calG(0) = \calG_0$ are equivalent.
\bigbreak

To define the notion of solitons in the generalized case, we first observe that, given an ECA~$(E, \pi, \ipr{\cdot}{\cdot} [\cdot, \cdot])$ and a scalar~$c \in \R^+$, the data
\begin{inline}
	(E, \pi, c \cdot \ipr{\cdot}{\cdot}, [\cdot, \cdot])
\end{inline}
defines an exact Courant algebroid again.
Following the notation in \cite{Fusi2024}, we will denote this ECA by $c \cdot E$.

\begin{definition}[{\cite[Def.~3.1]{Fusi2024}}] \label{def:grs}
A solution~$(\calG(t))_{t \in I}$ to the generalized Ricci flow~\labelcref{eq:grf} is called \emph{generalized Ricci soliton~(GRS)}, if there exists a one-parameter family of scalings~$c\colon I \to \R^+$ with $c(0) = 1$ together with a family of Courant algebroid isometries
\begin{inline}
	F_t\colon (c(t) \cdot E, \calG(0)) \overset{\cong}{\longto} (E, \calG(t)).
\end{inline}
\end{definition}

If $(E, \calG) \cong_\sigma (\mathbb{T}M_H, \calG(g, b))$ for some metric Courant algebroid~$(E, \calG)$ with isotropic splitting~$\sigma$, it follows that $\calG$ is a generalized metric on $c \cdot E$ as well, where
\begin{inline}
	(c \cdot E, \calG) \cong_\sigma (\mathbb{T}M_{cH}, \calG(c g, c b)).
\end{inline}
Therefore, equivalent to \cref{def:grs}, a family~$\calG(g(t), b(t))$ on ${\mathbb{T}M}_{H_0}$ is a generalized Ricci soliton if there exists a one-parameter family of diffeomorphisms~$\varphi_t$, such that
\begin{equation} \label{eq:grs_pull_back_rescaling}
	g(t) = c(t)\, \varphi_{-t}^* g(0)
	\qquad\text{and}\qquad
	H(t) = c(t)\, \varphi_{-t}^* H_0,
\end{equation}
where $H(t) \ceq H_0 + \D b(t)$.
Similar to the ordinary Ricci flow, we may derive the soliton solution to \cref{eq:grf_bismut} and obtain the equivalent infinitesimal conditions
\begin{equation} \label{eq:grs}
	\begin{cases}
		\Ric_{g_0, H_0}^B = \lambda g_0 + \frac{1}{2} L_X g_0, \\
		\Laplace^{g_0} H_0 = -2 \lambda H_0 - L_X H_0,
	\end{cases}
\end{equation}
on $\pa*{{\mathbb{T}M}_{H_0}, \calG(g_0, 0)}$.

\subsection{Algebraic Solitons on Lie Groups}

Over a Lie group~$\GG$, we wish to study Courant algebroid structures compatible with the group's left action.
As Lie groups are parallelizable, every ECA over~$\GG$ is isomorphic to
\begin{inline}
	\underline{\frakg \oplus \frakg^*} \ceq \GG \times \pa*{\frakg \oplus \frakg^*}
\end{inline}
as vector bundles, where $\frakg \ceq T_e \GG$ refers to the Lie algebra of $\GG$.
In the following, we abbreviate $\frakG \ceq \frakg \oplus \frakg^*$.

\begin{definition}
A metric Courant algebroid~$({\mathbb{T}\GG}_H, \calG(g, b))$ over a Lie group~$\GG$ is \emph{left-invariant} if the Riemannian metric~$g$, the two-form~$b \in \Omega^2(\GG)$ and the three-form~$H \in \Omega^3_{\mathrm{cl}}(\GG)$ are left-invariant.
\end{definition}
Equivalently, for all $g \in \GG$, the left action~$l_g\colon \GG \to \GG$ lifts to isometries
\begin{inline}
	\overline{l_g} \ceq \pa[\big]{\pa{l_g}_*, \pa[\big]{l_g^{-1}}^*}
\end{inline}
of the $H$-twisted generalized tangent bundle.
For the remainder of this paper, $({\mathbb{T}\GG}_H, \calG(g, b))$ will be assumed left-invariant, allowing us to restrict to the fiber over the identity in $\GG$.
For $[\cdot, \cdot]$ the Dorfman bracket on ${\mathbb{T}\GG}_H$, we use the notation $\frakG_{[\cdot, \cdot]}$ for the exact Courant algebroid~$(\frakG, \pi, [\cdot, \cdot], \ipr{\cdot}{\cdot})$ with $\pi(X + \xi) \ceq X$ and $\ipr{X + \xi}{X + \xi} \ceq \xi(X)$ for all $X + \xi \in \frakG$.
Note, that $(\frakG, [\cdot, \cdot], \ipr{\cdot}{\cdot})$ is a quadratic Lie algebra.
\medbreak

Following \cite{Fusi2024}, we consider the group
\begin{inline}
	\mathrm{L} \ceq \Lambda^2 \frakg^* \rtimes \GL(\frakg) \le_\iota \mathrm{O}(\frakG, \ipr{\cdot}{\cdot})
\end{inline}
which is embedded in the orthogonal maps via
\begin{equation*}
	\iota (b, h)
	= e^b \cdot
			\begin{bmatrix}
				h & 0 \\
				0 & (h^{-1})^*
			\end{bmatrix}.
\end{equation*}
Naturally, $\mathrm{L}$ acts on generalized metrics by $\ell \cdot \calG \ceq \ell \circ \calG \circ \ell^{-1}$ and on Dorfman brackets by
\begin{inline}
	\ell \cdot [\cdot, \cdot] \ceq \ell \circ \bra*{\ell^{-1} \cdot, \ell^{-1}}.
\end{inline}
It is immediately obvious from \cref{eq:gen_metric_TTM} that $\mathrm{L}$ acts transitively on the space of all left-invariant generalized metrics.
Thus, instead of fixing the Courant algebroid structure and flowing the generalized metric, on a Lie group, we may permanently fix a background metric~$\bar{\calG}$ and vary the Dorfman bracket instead.
To make this statement precise, a Dorfman bracket on $\frakG$ is a map~$\bmu \in \Lambda^2\frakG^* \tprod \frakG$ which satisfies
\begin{equation*}
    \ipr{\bmu(\cdot, \cdot)}{\cdot} \in \Lambda^3\frakg^*,\qquad
    \bmu|_{\frakg^* \tprod \frakg^*} = 0 \qquad\text{and}\qquad
    \mathcal{J}(\bmu) = 0
\end{equation*}
with $\mathcal{J}$ the Jacobiator on $\Lambda^2\frakG^* \tprod \frakG$.
We denote the set of all Dorfman brackets by $\calD \subset \Lambda^2\frakG^* \tprod \frakG$.

Furthermore, we denote the representation of the Lie algebra~$\mathfrak{l} = \operatorname{Lie}(\mathrm{L})$ on $\Lambda^2\frakG^* \tprod \frakG$ and of $\mathfrak{gl}(\frakg)$ on $\Lambda^2\frakg^* \tprod \frakg$ by
\begin{align*}
	\Theta\colon \mathfrak{l} &\longto \End\pa*{\Lambda^2\frakG^* \tprod \frakG},\quad
	\Theta(A) \bmu \ceq A\bmu(\cdot, \cdot) - \bmu(A \cdot, \cdot) - \bmu(\cdot, A\cdot), \\
	\theta\colon \mathfrak{gl}(\frakg) &\longto \End\pa*{\Lambda^2\frakg^* \tprod \frakg},\quad
	\theta(A) \mu \ceq A\mu(\cdot, \cdot) - \mu(A \cdot, \cdot) - \mu(\cdot, A\cdot)
\end{align*}
and the canonical representation on forms by
\begin{inline}
	\rho\colon \mathfrak{gl}(\frakg) \longto \End\pa*{\Lambda^*\frakg^*},
\end{inline}
where
\begin{equation*}
	\rho(A)\omega(\range{X_1}{X_k}) \ceq -\omega(\range{AX_1}{X_k}) - \dotsc -\omega(\range{X_1}{AX_k})
\end{equation*}
for all $\omega \in \Lambda^k\frakg^* \le \Lambda^*\frakg^*$ and all $\range{X_1}{X_k} \in \frakg$.

Now, referring to the Lie-bracket of $\frakg$ as $\mu$ and to the Dorfman bracket on $\frakG$ as $\bmu$, the latter then takes the shape
\begin{equation*}
	\bmu\pa*{X + \xi, Y + \eta}
	= \mu(X, Y) - \eta \circ \mu_X + \xi \circ \iota_Y + \iota_Y \iota_X H
\end{equation*}
with $\mu_X \ceq \mu(X, \cdot)$ and is thus uniquely determined by $H$ and $\mu$.
Conversely, we can recover $\mu$ and $H$ from the Dorfman bracket via the projections
\begin{inline}
	\mu = \bmu_{\frakg} \ceq \pi \circ \bmu
	\ \text{and}\
	H = \bmu_{\Lambda^3} \ceq 2\ipr{\bmu(\cdot, \cdot)}{\cdot}.
\end{inline}
Therefore, $\mathrm{L}$ acts on $\calD$ via
\begin{inline}
	\bra*{(b, h) \cdot \bmu}_{\frakg} = h \cdot \mu
	\ \text{and}\
	\bra*{(b, h) \cdot \bmu}_{\Lambda^3} = h \cdot (H - \D \beta)
\end{inline}
respectively.

All this yields the following flow:
\begin{definition}[{\cite[Def.~7.1]{Fusi2024}}] \label{def:l_norm_bf}
A family of Dorfman brackets~$(\bmu(t))_{t \in I}$ in $\Lambda^2\frakG^* \tprod \frakG$ on the background metric Courant algebroid~$(\frakG, \bar{\calG} \ceq \calG(\overline{g}, 0))$ solves the \emph{$\ell$-normalized generalized bracket flow} with initial value~$\bmu_0 \in \calD$, if
\begin{equation} \label{eq:l_norm_bf}
	\diffp*{\bmu}{t} = -\Theta\pa*{{\mathcal{R}c}_{\bmu} - A_{\bmu}} \bmu + \ell \bmu,\qquad \bmu(0) = \bmu_0,
\end{equation}
where $\ell = \ell(t)$ is a smooth function along the flow (depending on $t$), ${\mathcal{R}c}_{\bmu}$ denotes the generalized Ricci curvature and
\begin{inline}
	A_{\bmu} \ceq \operatorname{Skew}_{\overline{\calG}}\pa*{\D^* \bmu_{\Lambda^3}}.
\end{inline}
\end{definition}
\Cref{def:l_norm_bf} is a (slightly) more general definition than that of a (non-normalized) bracket flow, which we obtain by setting $\ell \equiv 0$, i.e.\
\begin{equation} \label{eq:ordinary_bf}
	\diffp*{\bmu}{t} = -\Theta\pa*{{\mathcal{R}c}_{\bmu} - A_{\bmu}} \bmu.
\end{equation}
Conversely, the unique solution~$\bmu^\ell$ to \cref{eq:l_norm_bf} and the (also unique) solution~$\bmu$ to \cref{eq:ordinary_bf} with same initial Dorfman bracket are related via rescaling in time and space, so $\bmu^\ell(t) = c(t) \bmu(\tau(t))$ for a scaling~$c$ and time reparametrization~$\tau$ along the flows.
Additionally, for a metric Courant algebroid~$(\frakG_{\bmu_0}, \calG_0)$, the bracket flow of $\bmu_0$ and the generalized Ricci flow of $\calG_0$ are equivalent in the following way:
If $(\calG(t))_{t \in [0, T)}$ is a family of generalized metrics solving \cref{eq:grf} and $(\bmu(t))_{t \in [0, \tilde{T})}$ a family solving \cref{eq:ordinary_bf} on their respective maximal time intervals of existence, then $T = \tilde{T}$ and there exists a family of generalized isometries
\begin{equation*}
	F_t\colon (\frakG_{\bmu_0}, \calG(t)) \overset{\cong}{\longto} (\frakG_{\bmu(t)}, \calG_0).
\end{equation*}
\medbreak

By extending the background metric~$\bar{\calG} = \calG(\bar{g}, 0)$ to tensors in $\Lambda^2 \frakG^* \tprod \frakG$ via
\begin{align*}
	\bar{g}(\bmu, \boldsymbol{\nu})
	&\ceq 2 \sum_{i, j}\ipr*{\bar{\calG}\bmu(e_i, e_j)}{\boldsymbol{\nu}(e_i, e_j)} \\
			&\qquad+ 4 \sum_{i, j}\ipr*{\bar{\calG}\bmu(e_i, e^j)}{\boldsymbol{\nu}(e_i, e^j)}
			+ 2 \sum_{i, j}\ipr*{\bar{\calG}\bmu(e^i, e^j)}{\boldsymbol{\nu}(e^i, e^j)}
\end{align*}
with $\set{e_i}$ and $\set{e^j}$ orthonormal bases of $\frakg$ and $\frakg^*$ respectively, we can choose the scaling~$c$ in such a way that $\bmu^\ell$ has constant norm along the flow.
This is achieved by
\begin{inline}
	c = \frac{1}{\abs{\bmu}},
\end{inline}
yielding
\begin{equation} \label{eq:ell_scalar_normalization}
	\ell = \ell_\bmu \ceq \frac{1}{\abs{\bmu}^4}\, \bar{g}\pa*{\Theta({\mathcal{R}c}_{\bmu} - A_\bmu) \bmu, \bmu}.
\end{equation}
We refer to the $\ell$-normalized bracket with normalization function \labelcref{eq:ell_scalar_normalization} as \emph{scalar-normalized bracket flow}, the fixed points of which are called \emph{algebraic solitons}:
\begin{definition}[{\cite[Def.~7.3]{Fusi2024}}]
Let $\bmu \in \calD$ be a Dorfman bracket on $\frakG$ with background metric~$\calG$.
Then $\bmu$ is said to be a \emph{generalized algebraic soliton} if it is a fixed point of the scalar-normalized bracket flow, that is,
\begin{inline}
	\Theta({\mathcal{R}c}_{\bmu} - A_\bmu) \bmu = \ell_\bmu \bmu
\end{inline}
for some $\ell_\bmu \in \R$.
The soliton is said to be \emph{expanding (steady, shrinking)}, if $\ell_\bmu > 0$ ($\ell_\bmu = 0$, $\ell_\bmu < 0$).
\end{definition}
Equivalent to this definition, $\bmu \in \calD$ is an algebraic soliton, if for $\mu = \bmu_{\frakg}$, $H = \bmu_{\Lambda^3}$,
\begin{equation}
	\label{eq:nil_soliton}
	\left\{
	\begin{aligned}
		\Ric^B_{\mu, H} &= \lambda \Id + D, \\
		\Laplace_\mu H &= \lambda H + \rho\pa[\Big]{\Ric^B_{\mu, H}} H,
	\end{aligned}
	\right.
\end{equation}
where $\lambda \in \R$ and $D \in \Der(\frakg)$ is a derivation of the Lie bracket~$\mu$.
Because of this, we will also say that the triple~$(g, H, \lambda)$ is a (generalized) algebraic soliton, if it solves \cref{eq:nil_soliton}.

Trivially, $(g, H, \lambda)$ is a solution to \cref{eq:nil_soliton} if and only if $(g, -H, \lambda)$ is.
Moreover, algebraic solitons can be rescaled in time and space, similar to ordinary Ricci solitons:
\begin{lemma} \label{lem:rescaling_soliton}
If $(g, H, \lambda)$ is an algebraic soliton, so is $\pa*{c g, c H, c^{-1} \lambda}$ for all $c \in \R^+$.
\end{lemma}
\begin{proof}
Let $\tilde{g} \ceq cg$, $\tilde{H} \ceq cH$.
Then $\tilde{g}^{-1} = c^{-1} g^{-1}$ and
\begin{equation*}
	\tilde{H}^2(X, Y)
	= \tilde{g}\pa{X \lrcorner \tilde{H}, Y \lrcorner \tilde{H}}
	= c^{-2} g \pa*{X \lrcorner cH, Y \lrcorner cH}
	= H^2(X, Y).
\end{equation*}
As furthermore $\Ric^{\tilde{g}} = \Ric^{cg} = \Ric^g$, we have
\begin{inline}
	\Ric^B_{\tilde{g}, \tilde{H}}
	= c^{-1} \Ric^B_{g, H}.
\end{inline}
For $\alpha, \beta \in \Lambda^k\frakg^*$, it is
\begin{equation*}
	\alpha \wedge \tilde{*} \beta
	= \tilde{g}(\alpha, \beta) \cdot \vol^{\tilde{g}}
	= c^{-k} g(\alpha, \beta) \cdot c^{\frac{n}{2}} \vol^g
	= c^{\frac{n}{2} - k} \pa*{\alpha \wedge \hodge \beta},
\end{equation*}
so
\begin{inline}
	\D^*_{\tilde{g}}
	= (-1)^s\, \tilde{*}\, {\D}\, \tilde{*}
	= c^{-1} \D^*_{g},%
\ \text{which yields}\
	\Laplace^{\tilde{g}} = -(\D \D^*_{\tilde{g}} + \D^*_{\tilde{g}} \D) = c^{-1} \Laplace^{g}.
\end{inline}
Thus, all quantities in \cref{eq:nil_soliton} scale in such a way that $\pa*{cg, cH, c^{-1} \lambda}$ is also a solution.
\end{proof}

\subsection{Correspondence Between Generalized Ricci Solitons and Algebraic Solitons}

We finish the preliminaries by stating a correspondence between left-invariant generalized Ricci solitons and (semi-)algebraic solitons.
\begin{definition} \label{def:semi_alg_soliton}
The triple~$(g, H, \lambda)$ is called \emph{semi-algebraic generalized soliton}, if
\begin{gather*}
	\left\{
	\begin{aligned}
		\Ric^B_{\mu, H} &= \lambda \Id + S(D), \\
		\Laplace_\mu H &= -2\lambda H + \rho\pa{D} H,
	\end{aligned}
	\right.
\end{gather*}
where $D \in \Der(\mu)$ is a derivation of the Lie-bracket and $S(D) = \frac{1}{2}(D + D^\intercal)$ is the symmetric part of $D$.
\end{definition}
An equivalent definition can also be found in \cite[Section~7.1]{Fusi2026}.
Clearly, all generalized algebraic solitons are also semi-algebraic, and semi-algebraic solitons are algebraic if $D$ is symmetric.
For ordinary Ricci solitons (i.e.\ $H = 0$), it is a well-known result that left-invariant Ricci solitons always correspond to algebraic solitons and semi-algebraic solitons are always algebraic \cite{Jablonski2014}.
In the following, we extend this result to the generalized case and show that similarly, generalized Ricci solitons, algebraic solitons and semi-algebraic solitons coincide on nilpotent Lie groups.
The crucial ingredients for this are augmented Dorfman brackets and the theory of gauged-flow and results therein of \cite[Section~4.2]{Fusi2026}.

\begin{lemma} \label{lem:nilpot_grs_is_semi-alg}
If $(g, H, \lambda)$ is a left-invariant generalized Ricci soliton on a nilpotent Lie group~$\NN$ with $\NN$ connected and simply-connected, $\frakn \ceq \operatorname{Lie}(\NN)$ the Lie algebra, it is also a semi-algebraic soliton.
\end{lemma}
\begin{proof}
We adapt the reasoning in \cite[Prop.~A.1]{He2015} to the generalized case.
By definition,
\begin{equation} \label{eq:grs2}
	\begin{cases}
		\Ric^g - \frac{1}{4} H_g^2 = \lambda g + \frac{1}{2} L_X g, \\
		-\D\D^*_g H = -2 \lambda H - L_X H,
	\end{cases}
\end{equation}
where $g$ and $H$ are left-invariant on $\NN$ and $X \in \VF(\NN)$ is the soliton vector field.
Consider the Lie group
\begin{gather*}
	S \ceq \operatorname{Iso}(\NN, g)_H = \iset{\varphi \in \operatorname{Iso}(\NN, g)}{\varphi^* H = H}
\shortintertext{with associated Lie algebra}
	\mathfrak{s} = \iset{Y \in \VF(\NN)}{L_Y g = 0,\, L_Y H = 0}.
\end{gather*}
For all $Y \in \mathfrak{s}$,
\begin{align*}
	L_{[X, Y]} g
	= L_X L_Y g - L_Y L_X g
	= 0 - 2 L_Y \pa*{\Ric^g - \frac{1}{4} H_g^2 - \lambda g} = 0
\shortintertext{and}
	L_{[X, Y]} H
	= L_X L_Y H - L_Y L_X H
	= 0 - L_Y \pa*{\D \D_g^* H - 2\lambda H} = 0,
\end{align*}
so $[X, Y] \in \mathfrak{s}$.
As a result, $\operatorname{ad}_X$ restricts to $\mathfrak{s}$, i.e.\ $\operatorname{ad}_X\colon \mathfrak{s} \to \mathfrak{s}$ is a derivation of $\mathfrak{s}$.
Let $\mathfrak{nil}(\mathfrak{s})$ denote the nilradical of $\mathfrak{s}$.
Because $\mathfrak{nil}(\mathfrak{s})$ is nilpotent, so is $\operatorname{ad}_X(\mathfrak{nil}(\mathfrak{s}))$, and by the maximality of $\mathfrak{nil}(\mathfrak{s})$, $\operatorname{ad}_X(\mathfrak{nil}(\mathfrak{s})) \subseteq \mathfrak{nil}(\mathfrak{s})$, which means that $\operatorname{ad}_X$ further restricts to $\mathfrak{nil}(\mathfrak{s})$.
On the other hand, due to \cite{Wilson1982}, we know that the Lie algebra~$\operatorname{Lie}(\NN^L) = \VF(\NN)^R$ of $\NN^L \ceq \iset{l_h}{h \in \NN} \subseteq \operatorname{Iso}(\NN, g)$ is exactly the nil-radical of $\operatorname{Lie}(\operatorname{Iso}(\NN, g))$, yielding
\begin{inline}
	\mathfrak{nil}(\mathfrak{s}) \subseteq \operatorname{Lie}(\NN^L).
\end{inline}
On the other hand, $g$ and $H$ are invariant under the left action of $\NN$, so
\begin{inline}
	 \operatorname{Lie}(\NN^L) \subseteq \mathfrak{nil}(\mathfrak{s})
\end{inline}
and thus $\mathfrak{nil}(\mathfrak{s})$ is the Lie algebra of right-invariant vector fields on $\NN$.
By setting $Z \ceq X - Y$ with $Y$ the right-invariant extension of $X|_e$, we obtain $L_Z g = L_X g$, $L_Z H = L_X H$ and
\begin{equation} \label{eq:ad_Z_Derivation_NL}
	\operatorname{ad}_Z \in \operatorname{Der}\pa*{\VF(\NN)^R} \qquad\text{with}\qquad Z|_e = 0.
\end{equation}
Let $\varphi^Z_t$ denote the flow of $Z$.
For every integral curve~$\sigma(t) \ceq \varphi^Z_t(p)$, we calculate the derivative of $\abs{Z}_{g}^2|_{\sigma(t)}$ in $t$,
obtaining
\begin{equation*}
	\diff{}{t} \abs{Z}^2
	= L_Z g(Z, Z)
	= 2 \pa*{\Ric^g(Z, Z) - \frac{1}{4} H_g^2(Z, Z) - \lambda \abs{Z}^2}.
\end{equation*}
As $\Ric^g$ and $H^2_g$ are left-invariant, they are bounded, allowing us to estimate
\begin{inline}
	\abs*{\partial_t \abs{Z}_g^2}
	\le C\, \abs{Z}_g^2
\end{inline}
for some constant~$C > 0$ along $\sigma(t)$.
As $(\NN, g)$ is complete, it thus follows analogously to the proof for ordinary Ricci solitons \cite[Thm.~2.25]{Chow2023} that the flow of $Z$ must be complete using Grönwall's inequality.
\smallbreak

By \cref{eq:ad_Z_Derivation_NL}, $\varphi^Z_{t}(e) \equiv e$ and $\varphi^Z_t$ normalizes $\NN^L$, so
\begin{equation*}
	\varphi^Z_{t} \circ l_h \circ \varphi^Z_{-t} = l_{\varphi^Z_{t}(h)}
\end{equation*}
and thus
\begin{equation*}
	\varphi^Z_{t}(a \cdot b) = \varphi^Z_{t}(a) \cdot \varphi^Z_{t}(b)
\end{equation*}
for all $a, b \in \NN$, i.e.\ $\varphi^Z_{t} \in \Aut(\NN)$.
As $\NN$ is simply-connected, $\Aut(\NN) \cong \Aut(\frakn)$, so we can find a derivation~$D \in \Der(\frakn)$ tangential to the curve~$t \mapsto \D\varphi^Z_t|_e \in \Aut(\frakn)$ in $t = 0$.
Now, on the Lie algebra,
\begin{align*}
	L_X g
	&= L_Z g
	= \left.\diffp{}{t}\right|_{t = 0} \pa*{\varphi^Z_t}^*g
	= g(D \cdot, \cdot) + g(\cdot, D\cdot)
	= 2g(S(D) \cdot, \cdot)
\shortintertext{and}
	L_X H
	&= L_Z H
	= \left.\diffp{}{t}\right|_{t = 0} \pa*{\varphi^Z_t}^*H
	= -\rho(D)H,
\end{align*}
where $S(D) \ceq \frac{1}{2}\pa*{D + D^\intercal}$.
By \cref{eq:grs}, we thus obtain
\begin{equation*}
	\Ric_{\mu, H}^B = \lambda \Id + S(D)
\qquad\text{and}\qquad
	\Laplace^{g} H = -2 \lambda H + \rho(D) H.
	\qedhere
\end{equation*}
\end{proof}

\begin{theorem}
\label{thm:nilpotent-geometric-algebraic}
Let $\GG$ be a connected, simply-connected Lie group, $g$ a left-invariant
Riemannian metric on $\GG$, and $H$ a left-invariant closed three-form.
\begin{enumerate}
\item Every algebraic generalized Ricci soliton $(g, H, \lambda)$ is a generalized Ricci soliton.
\item If $\GG=\mathrm{N}$ is nilpotent, the notions of generalized Ricci soliton, semi-algebraic and algebraic generalized Ricci soliton are equivalent.

Equivalently, there exists a $g$-symmetric derivation $D\in\operatorname{Der}(\frakn)$, $\mathfrak{n} \ceq \operatorname{Lie}(\mathrm{N})$ such that
\begin{equation}
\label{eq:nilpotent-algebraicity-characterization}
g^{-1} \Ric_{g, H}^B
=\lambda\operatorname{Id}+D,
\qquad
\D_g^*H=0,
\qquad
\rho(D)H=2\lambda H.
\end{equation}
In particular, the torsion of every left-invariant geometric generalized Ricci soliton is harmonic.
\end{enumerate}
\end{theorem}

\begin{proof}
Part~(i) is \cite[Lemma~7.5]{Fusi2024}.
For part~(ii), recall from \cref{eq:ricci_bismut} that
\begin{equation*}
	P \ceq g^{-1} \Ric_{g, H}^B = g^{-1}\left(\operatorname{Ric}^g-\frac14H_g^2\right)
\end{equation*}
is the symmetric Bismut Ricci operator as seen before, where the skew part of the full Bismut Ricci tensor is $-\frac12\D_g^*H$.
From \cref{lem:nilpot_grs_is_semi-alg}, we already know that all left-invariant generalized Ricci solitons on $\mathrm{N}$ are semi-algebraic.
In particular, we have
\begin{equation}
\label{eq:nilpotent-semi-algebraic-equations}
P=\lambda\operatorname{Id}+\frac12(A+A^\intercal),
\qquad
\Delta_gH=-2\lambda H+\rho(A)H, \quad A\in \operatorname{Der}(\frakn).
\end{equation}

We show that a semi-algebraic soliton is algebraic, by applying
\cite[Lem.\ 4.8, Thm.\ 4.12]{Fusi2026} to a suitable stationary
normalized flow.
Suppose that
\eqref{eq:nilpotent-semi-algebraic-equations} holds with
$A\in\operatorname{Der}(\frakn)$.
Put $\beta=\D_g^*H$.
Since $\D H=0$, its torsion equation becomes
\begin{equation}
\label{eq:nilpotent-torsion-defect}
\rho(A)H=2\lambda H-\D\beta.
\end{equation}

Suppose $\mu$ is the bracket of $\frakn$.
The Riemannian metric~$g$ on $\frakn$ extends to the metric~$h=g\oplus g^{-1}$ on $\frakN \ceq \frakn\oplus\frakn^*$, while the associated Dorfman
bracket on $\frakN$ is $\bmu$ given by
\begin{equation}
\label{eq:nilpotent-doubled-bracket}
\bmu(U+\xi,V+\eta)
=\mu(U,V)-\eta\circ\operatorname{ad}_U
+\xi\circ\operatorname{ad}_V+H(U,V,\cdot).
\end{equation}
Closedness of $H$ shows that $\bmu$ satisfies the Jacobi identity and hence is a Lie bracket.
Moreover, the Lie algebra~$(\frakN, \bmu)$ is nilpotent.
To see this, we apply Engel's theorem.
Let $U+\xi \in \frakN$.
Its adjoint operator has the form
\begin{equation*}
\operatorname{ad}^{\bmu}_{U+\xi} = \begin{pmatrix}
\operatorname{ad}^{\mu}_U & 0\\
\xi\circ \operatorname{ad}^{\mu}_V+H(U, V, \cdot) & -(\operatorname{ad}^{\mu}_U)^*
\end{pmatrix}.
\end{equation*}
As $\frakn$  is nilpotent, $\operatorname{ad}^{\mu}_U$ and its dual are nilpotent operators.
Thus, every adjoint map is block triangular with nilpotent diagonal blocks,
hence $\operatorname{ad}^{\bmu}_{U+\xi}$ is nilpotent and by Engel's theorem we get that $(\frakN, \bmu)$ is a nilpotent algebra.
In particular, its Killing form vanishes.
Now, recall from \cite[Section~4.1]{Fusi2026} that an augmented bracket on any $\frakG \ceq \mathfrak g\oplus \mathfrak g^*$ is a pair $(\bmu, C)$, $\bmu$ a bracket and $C \in \End(\frakG)$ an endomorphism such that $C^\intercal C=3\mathrm{B}_{\bmu}$ with $\mathrm B_{\bmu}$ the Killing form operator, $h(\mathrm B_{\bmu}Z, W)=\kappa_\bmu(Z, W)$ and $\kappa_\bmu$ the Killing form.
In our setting $\frakN$ is nilpotent, hence its Killing form vanishes and hence we get $C^\intercal C=0$.
Since the metric $h$ is positive-definite, for any $Z$, we get
\begin{equation*}
|CZ|^2_h=h(C^\intercal CZ, Z)=0 \quad \implies \quad C=0.
\end{equation*}
Thus, in the notation of \cite{Fusi2026}, the augmented bracket $(\bmu, C)$ in our case is just $(\bmu, 0)$.

If $\bmu=0$, then $\mu=0$ and $H=0$, so $D=-\lambda\operatorname{Id}$ gives
\eqref{eq:nilpotent-algebraicity-characterization}.
We may thus assume
$\bmu\neq0$.
The equations are homogeneous under
\begin{equation*}
(\mu,H,\lambda,A)\longmapsto
(r\mu,rH,r^2\lambda,r^2A),\qquad r>0,
\end{equation*}
so we may normalize $|\bmu|_h=1$, where $h=g\oplus g^{-1}$.
All tensor norms use full contraction.
We also use the neutral form
\begin{equation*}
b_0(U+\xi,V+\eta)=\xi(V)+\eta(U),
\end{equation*}
and also identify $\frakn$ and $\frakn^*$ using $g$ for the following
block calculations, and take transposes with respect to $h$.

We now want to use the gauged-flow formalism of \cite{Fusi2026}.
Recall for a metric Lie algebra $(\mathfrak g, \mu, g)$, the mean curvature vector $U_\mu$ is the unique vector in $\mathfrak g$ such that
\begin{equation*}
g(U_\mu, V) = \operatorname{tr}(\operatorname{ad}^\mu_V), \quad \text{for\ every\ } V\in \mathfrak g \iff U_\mu=g^{-1}\operatorname{tr} \operatorname{ad}_{(\cdot)}.
\end{equation*}
This vector measures the failure of the Lie algebra~$\mathfrak g$ to be unimodular.
Notice that if $\mathfrak g$ is nilpotent (as is the case with us), hence unimodular then $U_{\mu}=0$.

The following computation is done so as to use \cite[Lemma~4.8]{Fusi2026}.
Let $B$ be defined by $g(BU,V)=\beta(U,V)$, so $B^\intercal =-B$, and put
\begin{equation}
\label{eq:nilpotent-doubled-derivation}
\mathcal D
=\begin{pmatrix}
A&0\\
-B&-2\lambda\operatorname{Id}-A^\intercal
\end{pmatrix},
\qquad
\mathcal R
=\begin{pmatrix}
P&\frac12B\\
-\frac12B&-P
\end{pmatrix}.
\end{equation}
By \cite[Lemma~7.6]{Fusi2024}, equation
\eqref{eq:nilpotent-torsion-defect} gives
$\mathcal D\in\operatorname{Der}(\bmu)$.
We also have
\begin{equation*}
    P=\lambda\operatorname{Id}+\frac12(A+A^\intercal), \quad \mathcal D^\intercal =
\begin{pmatrix}
A^\intercal &B\\
0&-2\lambda\operatorname{Id}-A
\end{pmatrix}.
\end{equation*}
Thus, we get
\begin{equation*}
\operatorname{Sym}_h(\mathcal D)
=\frac12(\mathcal D+\mathcal D^\intercal)
=
\begin{pmatrix}
\frac12(A+A^\intercal)&\frac12B\\
-\frac12B&-2\lambda\operatorname{Id}-\frac12(A+A^\intercal)
\end{pmatrix}.
\end{equation*}
Adding $\lambda\operatorname{Id}_{\frakN}$ gives
\begin{align*}
\lambda\operatorname{Id}_{\frakN}
+\operatorname{Sym}_h(\mathcal D)
&=
\begin{pmatrix}
\lambda\operatorname{Id}+\frac12(A+A^\intercal)&\frac12B\\
-\frac12B&-\lambda\operatorname{Id}-\frac12(A+A^\intercal)
\end{pmatrix}\\
&=
\begin{pmatrix}
P&\frac12B\\
-\frac12B&-P
\end{pmatrix}
=\mathcal R.
\end{align*}
Thus,
\begin{equation}
\label{eq:nilpotent-doubled-soliton}
\mathcal R=\lambda\operatorname{Id}
+\operatorname{Sym}_h(\mathcal D).
\end{equation}
Since $\frakn$ is nilpotent hence unimodular, the curvature $\mathcal R$ is the same as the modified curvature $\mathcal{R}c^\star$ in \cite[14]{Fusi2026}.

Choose at $(\bmu,0)$ the gauge value
\begin{equation}
\label{eq:nilpotent-gauge-value}
Q_0=-\operatorname{Skew}_h(\mathcal D)
=\begin{pmatrix}
-\frac12(A-A^\intercal)&\frac12B\\
\frac12B&-\frac12(A-A^\intercal)
\end{pmatrix}.
\end{equation}
This is an admissible gauge value in \cite[Lemma~4.8]{Fusi2026} and its free
orthogonal term is $-\frac12(A-A^\intercal)$.
More explicitly, for a general augmented
bracket $(\boldsymbol{\eta},Z)$, choose the free term in \cite[Lemma~4.8]{Fusi2026} to be
\begin{equation*}
q_{(\boldsymbol{\eta},Z)}=-\frac 12\pa*{|\boldsymbol{\eta}|_h^2+|Z|_h^2}(A-A^\intercal).
\end{equation*}
The lemma then gives a smooth gauge map $Q$, homogeneous of degree two,
whose value at $(\bmu,0)$ is $Q_0$.

Equations \eqref{eq:nilpotent-doubled-soliton} and \eqref{eq:nilpotent-gauge-value} imply
\begin{inline}
\mathcal R-Q_0=\lambda\operatorname{Id}+\mathcal D.
\end{inline}
Recall that the infinitesimal action of an endomorphism $T$ on the bracket $\bmu$ is given  by $\Theta$ which is described as
\begin{equation*}
\Theta(T)\bmu=T\bmu-\bmu(T\cdot,\cdot)-\bmu(\cdot,T\cdot).
\end{equation*}
The quantity $\Theta(T)\bmu$ measures the failure of the endomorphism $T$ to be a derivation.
Since $\mathcal D$ is a derivation, we have $\Theta(\mathcal D)\bmu=0$ and $\Theta(\operatorname{Id})\bmu=-\bmu$.
Consequently,
\begin{equation}
\label{eq:nilpotent-fixed-gauge}
\Theta(\mathcal R-Q_0)\bmu=-\lambda\bmu.
\end{equation}
Because $Q_0$ is $h$-skew-symmetric, its action on bracket tensors is
skew-symmetric.
Thus the normalizing coefficient in
\cite[Definition~4.9]{Fusi2026} is
\begin{inline}
\ell=\langle\Theta(\mathcal R)\bmu,\bmu\rangle_h=-\lambda,
\end{inline}
and the normalized gauged flow satisfies
\begin{equation*}
-\Theta(\mathcal R-Q_0)\bmu+\ell\bmu=0.
\end{equation*}
Its augmentation remains zero.
Hence $(\bmu,0)$ is a fixed point.
In fact, this is why we introduced the gauge in first place so that our bracket $\bmu$ becomes a fixed point of the normalized flow when viewed as an augmented bracket $(\bmu, 0)$.

The associated non-normalized gauged flow is immortal, i.e., it exists for all positive times by \cite[Theorem~A]{Fusi2024} and \cite[Theorem~4.6]{Fusi2026}, and its
normalized bracket is constantly $(\bmu,0)$; here we use the degree-two
homogeneity of the gauge and \cite[Remark~4.10]{Fusi2026}.
Therefore the limiting conclusion of \cite[Theorem~4.12]{Fusi2026} which describes normalized sub-sequential limits of non-zero immortal gauged flows,
applies to $(\bmu,0)$ itself.
The conclusion in \cite[Theorem~4.12]{Fusi2026}, in our notations, is that
\begin{equation*}
\mathcal R_{{\bmu_{\infty}}} + \ell_{\infty}\operatorname{Id}\in \operatorname{Der}(\bmu_{\infty}) \quad \text{and}\quad \D^*_{\mu_{\infty}}H_{\infty} = U_{\mu_{\infty}}\lrcorner H_{\infty}, \qquad \ell_{\infty}>0.
\end{equation*}
Here $\mathcal R$ is the generalized unimodular curvature, $\ell_{\infty}, \bmu_{\infty}, \mu_{\infty}, H_{\infty}$ are the limiting scale, bracket on $\frakN$ with corresponding bracket on $\frakn$ and the three-form respectively which are obtained after applying \cite[Theorem~4.12]{Fusi2026} and $U_{\mu_{\infty}}$ is the mean curvature vector.

Since we proved above that $(\bmu, 0)$ is an admissible augmented Dorfman bracket and the flow is constant, the limiting solution is $\bmu_{\infty}=\bmu$ itself.
Due to nilpotency, the mean-curvature vector~$U_{\mu}$ is zero, and hence we get
\begin{equation*}
\D_g^*H=0,\qquad
\mathcal R+\ell_\infty\operatorname{Id}
\in\operatorname{Der}(\bmu)
\end{equation*}
for some $\ell_\infty>0$.
Pairing $\Theta(\mathcal R)\bmu=\ell_\infty\bmu$ with the unit bracket
$\bmu$ identifies $\ell_\infty=-\lambda$.
Thus
\begin{equation}
\label{eq:nilpotent-symmetric-doubled-derivation}
\mathcal R-\lambda\operatorname{Id}
\in\operatorname{Der}(\bmu).
\end{equation}

Recall that we had $g(BU,V)=\beta(U,V)=\D^*H(U,V)$, hence $B=0$.
Clearly $P-\lambda\operatorname{Id}$ is
$g$-symmetric.
The derivation in
\eqref{eq:nilpotent-symmetric-doubled-derivation}, using \eqref{eq:nilpotent-doubled-derivation} has the form
\begin{equation*}
\begin{pmatrix}
P-\lambda\operatorname{Id} &0\\
0& -\lambda \Id-P
\end{pmatrix}.
\end{equation*}
The block characterization in \cite[Lemma~7.6]{Fusi2024} therefore yields
\begin{equation*}
P-\lambda\operatorname{Id}\in\operatorname{Der}(\frakn),
\qquad
\rho(D)H=2\lambda H.
\end{equation*}
Thus, $A^\intercal = 2(P-\lambda\operatorname{Id})-A\in \operatorname{Der}(\frakn)$.
This proves \eqref{eq:nilpotent-algebraicity-characterization} for every
semi-algebraic soliton.
These conditions imply the algebraic soliton equations, since
\begin{equation*}
\Delta_gH=0
=-2\lambda H+\rho(D)H
=\lambda H+\rho(P)H.
\end{equation*}
Conversely, an algebraic soliton is semi-algebraic by its defining
equations and $\rho(\operatorname{Id})H=-3H$.
Together with part~(i) and the geometric reduction above, this proves
all assertions.
\end{proof}

\begin{remark}
As explained throughout the proof, the nilpotent algebraic implication above is a fixed-point consequence of the gauged-flow results in \cite{Fusi2026}.
\end{remark}

\section{Generalized Nilsolitons on the Five-Dimensional Heisenberg Algebra} \label{sec:heisenberg_algebras}

\begin{definition} \label{def:heisenberg_algebra}
Let $V$ be a finite-dimensional real vector space equipped with a non-degenerate two-form~$\omega$.
The \emph{Heisenberg algebra}~$\heis(V, \omega)$ over the vector space~${V \oplus \R_z}$ is given by the bracket
\begin{equation*}
	\mu(X, Y)
	\ceq [X, Y]
	\ceq \omega(X, Y) \cdot z.
\end{equation*}
for all $X, Y \in \heis(V, \omega)$.
Furthermore, $k \ceq \frac{1}{2} \dim V$ is referred to as the \emph{rank} of $\heis(V, \omega)$.
\end{definition}
In the following, we fix a standard symplectic basis~$\set{p_i, q_i}_{i = 1}^n$ of $V$.
We will also use the notation~$\heis_{2k + 1} \ceq \heis\pa{\R^{2k}, \sum_i p^i \wedge q^i}$, where
\begin{equation} \label{eq:ordered_basis_heis}
	e_i
	\ceq
			\begin{cases}
				p_i, & 1 \le i \le k \\
				q_{i - k}, & k < i \le 2k \\
				z, & i = 2k + 1
			\end{cases}
\end{equation}
will be the ordered standard basis which we will work in throughout this section.
The Lie bracket is thus given by $[p_i, q_i] = z$ for all $i \in \setrange{1}{k}$ and zero otherwise.

Equivalent to \cref{def:heisenberg_algebra}, the Heisenberg algebra is given by
\begin{inline}
	\heis_n = \pa*{0, \cdots, 0, -\omega}
\end{inline}
in Salamon notation.

The closed forms on $\heis_n$ are very easy to compute:
\begin{lemma} \label{lem:H_closed}
Let $\heis_n = V \oplus \R_z$ denote the Heisenberg algebra of dimension~$n$ and let $k$ denote its rank (i.e.\ $n = 2k + 1$), with $z \in \heis_n$ a generator of the center~$\R_z$ and $\omega \ceq \sum_{i = 1}^k p^i \wedge q^i$ the symplectic form on $V$.
Additionally, let $L \ceq \omega \wedge \cdot$ denote the Lefshetz operator on $V$.

Then the Chevalley-Eilenberg differential on $\heis_n$ is given by
\begin{inline}
	\D = - L \circ (z\lrcorner).
\end{inline}
In particular, $\alpha \in \Lambda^\ell\heis_n^*$ is closed if and only if $\omega \wedge (z \lrcorner \alpha) = 0$.
Furthermore, for degree~$\ell \le k$, the closed forms are given by $\Lambda_{\mathrm{cl}}^\ell\heis_n^* = \Lambda^\ell V^*$.
\end{lemma}
\bigbreak

In \cite{Fusi2024}, it was already shown that the unique, non-classical nilsoliton~$(g, H, \lambda)$ on the three-dimensional Heisenberg algebra up to rescaling and automorphisms is the triple
\begin{equation*}
    g = p^* \tprod p^* + q^* \tprod q^* + z^* \tprod z^*,\qquad H = p^* \wedge q^* \wedge z^*,\qquad \lambda = -2.
\end{equation*}
We wish to extend this result to $\heis_5$.
As mentioned in the introduction, the classification will be done the in following way:
\begin{enumerate*}
\item In the first step, we use the automorphism group to simplify the metric and use this to
\item prove all solitons on $\heis_5$ necessarily have harmonic torsion.
This will then allow us to \item reason that after rescaling, for any nilsoliton, the standard basis defined above is orthonormal.
\item Next, we use the remaining freedoms given by the stabilizer~$\Aut(\heis_5)_g$ of this metric, to simplify the three-form~$H$, by complexifying $V = \R^4$ and studying the actions of the subgroups~$\SU(2)$ and $\UU(1)$ on the respective $(p, q)$-forms in $\Lambda^2_\C V^*$.
This reduces the soliton equations to two polynomials of degree three in two variables.
Solving these allows us to finally name all solutions.
\end{enumerate*}

\subsection{Canonical Form of the Metric}

We now classify all positive definite, symmetric bilinear forms~$g$ on Heisenberg algebras.
In the following lemma, we first describe the (well-known) automorphism group of the Heisenberg algebra:
\begin{lemma} \label{lem:aut_heisk}
The Lie algebra automorphisms of the $n = 2k + 1$-dimensional Heisenberg algebra are given by
\begin{equation*}
	\Aut(\heis_n)
	= \operatorname{Aff}_{\operatorname{CSp}}(\R^{2k})
	= \R^{2k} \rtimes \operatorname{CSp}(2k, \R),
\end{equation*}
where
\begin{inline}
	\operatorname{CSp}(2k, \R)
	\ceq \iset{A \in \GL(\R^{2k})}{\exists c \in \R^\times\colon A^* \omega = c \cdot \omega}
\end{inline}
denotes the (linear) conformal symplectic group.
In other words, the automorphism group of the Heisenberg algebra is the group of affine conformal symplectic transformations of $\R^{2k}$.

More explicitly, under the canonical identification $\heis_n = \R^{n}$ as vector spaces, the automorphism group~$\Aut(\heis_n) \le \GL(n, \R)$ has two connected components,
\begin{equation*}
	\Aut(\heis_n) = \Aut(\heis_n)^0 \sqcup \Aut(\heis_n)^0 \bar{\Omega},
\quad\text{where}\quad
	\bar{\Omega}
	\ceq
			\begin{bsmallmatrix}
				0 & I_k & 0 \\
				I_k & 0 & 0 \\
				0 & 0 & -1
			\end{bsmallmatrix}
\end{equation*}
and every element
\begin{equation*}
	A
	\in \Aut(\heis_n)^0
	= \operatorname{Aff}_{\operatorname{CSp}^+}(\R^{2k})
	= \iset{\varphi \in \operatorname{Aff}(\R^{2k})}{\varphi^* \omega = c \cdot \omega,\, c \in \R^+}
\end{equation*}
of the identity component uniquely decomposes into
\begin{equation} \label{eq:decomp_auto}
	A
	=
			\begin{bmatrix}
				I_{2k} & 0 \\
				d^\intercal & 1
			\end{bmatrix}
			\cdot
			\begin{bmatrix}
				S & 0 \\
				0 & 1
			\end{bmatrix}
			\cdot
			\begin{bmatrix}
				\sqrt{c} \cdot I_{2k} & 0 \\
				0 & c
			\end{bmatrix}
	\quad\text{with}\ d \in \R^{2k},\, S \in \Sp(2k, \R),\, c \in \R^+.
\end{equation}
\end{lemma}

From \cref{eq:decomp_auto}, we can directly see that a general derivation~$a \in \Der(\heis_n)$ is of the shape
\begin{equation} \label{eq:decomp_der}
	a
	=
			\begin{bmatrix}
				s + \frac{1}{2} c\, I_{2k} & 0 \\
				d^\intercal & c
			\end{bmatrix}
	\quad\text{with}\ d \in \R^{2k},\, s \in \mathfrak{sp}(2k, \R),\, c \in \R.
\end{equation}
\bigbreak

Using symplectic diagonalization, we obtain a very easy normal form of metrics on Heisenberg algebras:
\begin{lemma} \label{lem:simplification_metric}
Given an arbitrary, positive definite metric~$g \in \Sym_+^2 \heis_n^*$, there exists an automorphism~$A \in \Aut(\heis_n)$, such that in the ordered basis~\eqref{eq:ordered_basis_heis}, the coefficients are given by
\begin{inline}
	[(A^*g)\indices{_{ij}}]
	=
		\Diag(D, D, 1),
\end{inline}
where $D \in \R^{k \times k}$ with $n = 2k + 1$ is a positive definite diagonal matrix.
Furthermore, up to reordering, the matrix~$D$ is uniquely determined by $g$.
\end{lemma}
\begin{proof}
Let
\begin{equation*}
	[\tensor{g}{_{ij}}]
	=
		\begin{bmatrix}
			G_1 & c \\
			c^\intercal & g_2
		\end{bmatrix},\quad
\text{where}\
	G_1 \in \R^{2k \times 2k},\, c \in \R^{2k},\, g_2 \in \R^+.
\end{equation*}
let $S$, $d$ and $c$ be defined as in \cref{eq:decomp_auto}.

By the action of $\Aut(\heis_n)^0 \ge \R^+ \curvearrowright \Sym_+^2 \heis_n^*$, we may assume $g_2 = 1$.
Next, choosing $d \ceq -c \in \R^{2k} \le \Aut(\heis_n)$, the action on $g$ is
\begin{equation} \label{eq:conjugation_metric}
	\begin{bmatrix}
		I_{2k} & 0 \\
		d^\intercal & 1
	\end{bmatrix}^\intercal
			\cdot
			\begin{bmatrix}
				G_1 & c \\
				c^\intercal & 1
			\end{bmatrix}
			\cdot
			\begin{bmatrix}
				I_{2k} & 0 \\
				d^\intercal & 1
			\end{bmatrix}
	=
			\begin{bmatrix}
				G_1 - c c^\intercal & 0 \\
				0 & 1
			\end{bmatrix}.
\end{equation}
From \cref{eq:conjugation_metric}, it is also apparent that $G_1 - c c^\intercal$ is symmetric and positive definite.
It is a well-known result that every positive definite symmetric matrix can be diagonalized by a symplectic matrix, where the diagonal matrix is of the shape~$\Diag(D, D)$ and is unique up to reordering of the diagonal entries of $D$.
This is commonly referred to as \emph{Williamson's Theorem} and can e.g.\ be found in \cite{Williamson1936}.
Choosing $S \in \Sp(2k, \R)$ such that
\begin{equation*}
	S^\intercal \pa*{G_1 - c c^\intercal} S
	= \Diag(D, D)
\end{equation*}
with $D$ a positive definite diagonal matrix proves the theorem.
\end{proof}

As a result, we may without loss of generality assume
\begin{inline}
	[\tensor{g}{_{ij}}] = \Diag\pa*{D, D, 1}
\end{inline}
in the standard base.

\subsection{Classification of Nilsolitons on the Five-Dimensional Heisenberg Algebra}

The normal form of metrics on Heisenberg algebras from \cref{lem:simplification_metric} can be applied to show the following:
\begin{lemma} \label{lem:heis5_solitons_harmonic}
All nilsolitons on the five-dimensional Heisenberg algebra have harmonic torsion.
\end{lemma}
\begin{proof}
Let $(g, H, \lambda)$ denote a solution to \cref{eq:nil_soliton}.
Our aim is to prove $\Delta^g H = 0$.

By \cref{lem:simplification_metric}, we may w.l.o.g.\ assume that in the canonical basis of $\heis_5$, the metric is of the shape
\begin{inline}
	[\tensor{g}{_{ij}}] = \Diag\pa*{\lambda_1, \lambda_2, \lambda_1, \lambda_2, 1}
\end{inline}
with $\lambda_1, \lambda_2 > 0$.

The three-form~$H$ decomposes into $H = \beta \wedge z^* + H_0$, where $\beta \in \Lambda^2 V^*$ and $H_0 \in \Lambda^3 V^*$.
Furthermore, $(V, \omega)$ is a four-dimensional symplectic space, so we may further simplify $H_0 = \theta \wedge \omega$ with $\theta \in V^*$.

From the shape of the matrix~$[\tensor{g}{_{ij}}]$, we observe that $\UU(1) \times \UU(1) \le \Aut(\heis_5)_g$, where both $\UU(1)$-groups act via rotation in the respective $p$-$q$-planes.
Every element~$(\phi, \psi) \in \UU(1) \times \UU(1)$ fixes $z^*$ and $\omega$, so the action $\UU(1) \times \UU(1) \curvearrowright \Lambda^3 \heis_5^*$ also preserves the decomposition of $H$ above, in the sense that
\begin{inline}
	(\phi, \psi)^* H
	= \pa*{(\phi, \psi)^* \beta} \wedge z^* + \pa*{(\phi, \psi)^* \theta} \wedge \omega.
\end{inline}
Hence, the one-form~$\theta$ can be written as $\theta = \theta_1 p^1 + \theta_2 p^2$ with coefficients~$\theta_1, \theta_2 \ge 0$.
From this, we calculate
\begin{equation*}
	\Delta^g H
	= -\pa*{\frac{\theta_1}{\lambda_2^2}\, p^1 + \frac{\theta_2}{\lambda_1^2}\, p^2} \wedge \omega,
\end{equation*}
so the claim will follow if we are able to prove $\theta_1 = \theta_2 = 0$.

For the sake of a contradiction, assume this is not the case.
We start with the special case of neither $\theta_1$ nor $\theta_2$ vanishing.
The first condition of \cref{eq:nil_soliton} is equivalent to $\Ric^{1, 1} - \frac{1}{4} g^{-1} H^2 - \lambda \Id_5$ being of the shape of \cref{eq:decomp_der}.
Matching components in
\begin{equation*}
	\bra*{\tensor{\Ric}{^i_j} - \frac{1}{4} \tensor{g}{^{ik}} \pa*{\beta \wedge z^* + \pa*{\theta_1 p^1 + \theta_2 p^2} \wedge \omega}\indices{^2_{kj}} - \lambda \tensor{\delta}{^i_j}} \\
	=
			\begin{bmatrix}
				s + \frac{1}{2} c\, I_{2k} & 0 \\
				d^\intercal & c
			\end{bmatrix}
	=\vcentcolon a
\end{equation*}
with $d \in \R^{2k},\, s \in \mathfrak{sp}(2k, \R),\, c \in \R$, we obtain
\begin{align}
	c &= -{\lambda} - \frac{1}{4} \norm{\beta}_g^2 + \frac{{\lambda_1}^{2} + {\lambda_2}^{2}}{2 \, {\lambda_1}^{2} {\lambda_2}^{2}}, \notag \\
	d^\intercal &=
			\begin{bmatrix}
				\frac{{\beta_{13}} {\theta_1} {\lambda_1} + {\beta_{23}} {\theta_2} {\lambda_2}}{2 \, {\lambda_1} {\lambda_2}^{2}} &
				\frac{{\beta_{14}} {\theta_1} {\lambda_1} - {\beta_{13}} {\theta_2} {\lambda_2}}{2 \, {\lambda_1}^{2} {\lambda_2}} &
				\frac{{\beta_{12}} {\theta_2}}{2 \, {\lambda_1} {\lambda_2}} &
				-\frac{{\beta_{12}} {\theta_1}}{2 \, {\lambda_1} {\lambda_2}}
			\end{bmatrix} \label{eq:derivation_d_comp} \\
	0 &=
			\begin{bmatrix}
				\frac{{\beta_{13}} {\theta_1}}{2 \, {\lambda_1} {\lambda_2}^{2}} + \frac{{\beta_{23}} {\theta_2}}{2 \, {\lambda_1}^{2} {\lambda_2}} \\
				\frac{{\beta_{14}} {\theta_1}}{2 \, {\lambda_1} {\lambda_2}^{2}} - \frac{{\beta_{13}} {\theta_2}}{2 \, {\lambda_1}^{2} {\lambda_2}} \\
				\frac{{\beta_{12}} {\theta_2}}{2 \, {\lambda_1}^{2} {\lambda_2}} \\
				-\frac{{\beta_{12}} {\theta_1}}{2 \, {\lambda_1} {\lambda_2}^{2}}
			\end{bmatrix}, \label{eq:zero_obstr} \\
	s &=
		\frac{1}{2 \lambda_1^2 \lambda_2^2}
			\begin{bNiceMatrix}
				A & B \\
				C & D
			\end{bNiceMatrix}, \notag
\end{align}
where $A, B, C$ and $D$ are the block matrices
\begin{align*}
	A_{\bullet1} &=
		\begin{bsmallmatrix}
			-{\theta_1}^{2} {\lambda_1} - \theta_2^2 \lambda_2 + \frac{1}{2} \, {\left({\beta_{13}}^{2} - 1\right)} {\lambda_1}^{2}  - \frac{1}{2} \, {\left(2 \, {\lambda} {\lambda_1}^{2} + {\beta_{13}}^{2} + 3\right)} {\lambda_2}^{2}  - \frac{1}{2} \, {\left({\beta_{12}}^{2} + {\beta_{14}}^{2} - {\beta_{23}}^{2} - {\beta_{34}}^{2}\right)} {\lambda_1} {\lambda_2} \\
			{\beta_{13}} {\beta_{14}} {\lambda_1}^{2} - {\beta_{13}} {\beta_{23}} {\lambda_1} {\lambda_2}
		\end{bsmallmatrix}, \\
	A_{\bullet2} &=
		\begin{bsmallmatrix}
			{\beta_{13}} {\beta_{14}} {\lambda_1} {\lambda_2} - {\beta_{13}} {\beta_{23}} {\lambda_2}^{2} \\
			-{\theta_1}^{2} {\lambda_1} - {\theta_2}^{2} {\lambda_2} - \frac{1}{2} \, {\left({\beta_{13}}^{2} + 3\right)} {\lambda_1}^{2}  - \frac{1}{2} \, {\left(2 \, {\lambda} {\lambda_1}^{2} - {\beta_{13}}^{2} + 1\right)} {\lambda_2}^{2} - \frac{1}{2} \, {\left({\beta_{12}}^{2} - {\beta_{14}}^{2} + {\beta_{23}}^{2} - {\beta_{34}}^{2}\right)} {\lambda_1} {\lambda_2}
		\end{bsmallmatrix}, \\
	B &=
			\begin{bsmallmatrix}
			{\left({\beta_{12}} {\beta_{23}} - {\beta_{14}} {\beta_{34}}\right)} {\lambda_1} {\lambda_2} & -{\beta_{12}} {\beta_{13}} {\lambda_1} {\lambda_2} + {\beta_{13}} {\beta_{34}} {\lambda_2}^{2} \\
			{\beta_{13}} {\beta_{34}} {\lambda_1}^{2} - {\beta_{12}} {\beta_{13}} {\lambda_1} {\lambda_2} & -{\left({\beta_{12}} {\beta_{14}} - {\beta_{23}} {\beta_{34}}\right)} {\lambda_1} {\lambda_2}
			\end{bsmallmatrix}, \\
	C &=
			\begin{bsmallmatrix}
			{\left({\beta_{12}} {\beta_{23}} - {\beta_{14}} {\beta_{34}}\right)} {\lambda_1} {\lambda_2} & {\beta_{13}} {\beta_{34}} {\lambda_1} {\lambda_2} - {\beta_{12}} {\beta_{13}} {\lambda_2}^{2} \\
			-{\beta_{12}} {\beta_{13}} {\lambda_1}^{2} + {\beta_{13}} {\beta_{34}} {\lambda_1} {\lambda_2} & -{\left({\beta_{12}} {\beta_{14}} - {\beta_{23}} {\beta_{34}}\right)} {\lambda_1} {\lambda_2}
			\end{bsmallmatrix},
\allowdisplaybreaks[1]
\shortintertext{and}
	D_{\bullet1} &=
		\begin{bsmallmatrix}
			\frac{1}{2} \, {\left({\beta_{13}}^{2} - 1\right)} {\lambda_1}^{2} - \theta_2^2 \lambda_2  - \frac{1}{2} \, {\left(2 \, {\lambda} {\lambda_1}^{2} + {\beta_{13}}^{2} + 3\right)} {\lambda_2}^{2} + \frac{1}{2} \, {\left({\beta_{12}}^{2} + {\beta_{14}}^{2} - {\beta_{23}}^{2} - {\beta_{34}}^{2}\right)} {\lambda_1} {\lambda_2} \\
			{\beta_{13}} {\beta_{23}} {\lambda_1}^{2} - {\beta_{13}} {\beta_{14}} {\lambda_1} {\lambda_2} + {\theta_1} {\theta_2} {\lambda_1}
		\end{bsmallmatrix}, \\
	D_{\bullet2} &=
		\begin{bsmallmatrix}
			-{\beta_{13}} {\beta_{14}} {\lambda_2}^{2} + {\left({\beta_{13}} {\beta_{23}} {\lambda_1} + {\theta_1} {\theta_2}\right)} {\lambda_2} \\
			-{\theta_1}^{2} {\lambda_1} - \frac{1}{2} \, {\left({\beta_{13}}^{2} + 3\right)} {\lambda_1}^{2}  + \frac{1}{2} \, {\left({\beta_{12}}^{2} - {\beta_{14}}^{2} + {\beta_{23}}^{2} - {\beta_{34}}^{2}\right)} {\lambda_1} {\lambda_2}  - \frac{1}{2} \, {\left(2 \, {\lambda} {\lambda_1}^{2} - {\beta_{13}}^{2} + 1\right)} {\lambda_2}^{2}
		\end{bsmallmatrix}.
\end{align*}
In order for $s$ to be in the Lie algebra~$\mathfrak{sp}(4, \R)$, necessarily $s^\intercal J + J s = 0$ for the canonical almost complex structure
\begin{equation*}
	J \ceq
			\begin{bmatrix}
				0 & -I_2 \\
				I_2 & 0
			\end{bmatrix},
\end{equation*}
which translates to $C = C^\intercal$, $B = B^\intercal$ and $A = -D^\intercal$.
From $A = -D^\intercal$, we get
\begin{align} \label{eq:A_-DT}
	-{\beta_{13}} {\beta_{14}} {\lambda_1} {\lambda_2} + {\beta_{13}} {\beta_{23}} {\lambda_2}^{2} = {\beta_{13}} {\beta_{23}} {\lambda_1}^{2} - {\beta_{13}} {\beta_{14}} {\lambda_1} {\lambda_2} + {\theta_1} {\theta_2} {\lambda_1}
\shortintertext{and}
	-{\beta_{13}} {\beta_{14}} {\lambda_1}^{2} + {\beta_{13}} {\beta_{23}} {\lambda_1} {\lambda_2} = -{\beta_{13}} {\beta_{14}} {\lambda_2}^{2} + {\left({\beta_{13}} {\beta_{23}} {\lambda_1} + {\theta_1} {\theta_2}\right)} {\lambda_2}
	\notag
\end{align}
as necessary conditions.
By assumption, neither $\theta_1$ nor $\theta_2$ vanish, so \cref{eq:A_-DT} shows $\beta_{13} \ne 0$ and $\lambda_1 \ne \lambda_2$.
Both together show
\begin{equation} \label{eq:obstr_theta1_theta2}
	\beta_{14} \lambda_1 - \beta_{23} \lambda_2 = 0.
\end{equation}
From \cref{eq:zero_obstr}, we know
\begin{equation*}
	\left\{\begin{aligned}
		\beta_{14} &= \beta_{13} \cdot \frac{\theta_2}{\theta_1} \cdot \frac{\lambda_2}{\lambda_1}, \\
		\beta_{23} &= -\beta_{13} \cdot \frac{\theta_1}{\theta_2} \cdot \frac{\lambda_1}{\lambda_2},
	\end{aligned}\right.
\end{equation*}
which we may plug into \cref{eq:obstr_theta1_theta2} and obtain
\begin{equation*}
	0
	= \frac{\theta_2}{\theta_1} \lambda_2 + \frac{\theta_1}{\theta_2} \lambda_1.
\end{equation*}
This contradicts the assumption $\theta_1 > 0$ and $\theta_2 > 0$, so at least one of those coefficients has to vanish.
By optionally applying an automorphism of the Heisenberg algebra that interchanges $p^1 \leftrightarrow p^2$, $q^1 \leftrightarrow q^2$ and $\lambda_1 \leftrightarrow \lambda_2$, we can w.l.o.g.\ assume $\theta_2 = 0$ and $\theta_1 > 0$.
Accordingly, \cref{eq:zero_obstr} implies
\begin{equation*}
	\beta_{12} = \beta_{13} = \beta_{14} = 0
\end{equation*}
and by \cref{eq:derivation_d_comp}, this also means $d = 0$, so the derivation~$a$ preserves the splitting~$V \oplus \R_z$.

Consequently, from the second condition of \cref{eq:nil_soliton}, we calculate
\begin{align*}
	0
	&= -\Laplace H + \lambda H + \rho\pa[\Big]{\Ric_{\mu, H}^B} H \\
	&= -\Laplace H - 2\lambda H + \rho\pa[\Big]{\Ric_{\mu, H}^B - \lambda \Id} H \\
	&= \pa*{\frac{1}{\lambda_2^2} - 2\lambda} \theta_1\, p^1 \wedge \omega
			- 2\lambda\, \beta \wedge z^*
			+ \rho\pa{a} \pa*{\beta \wedge z^*} + \theta_1\, \rho(a) \pa{p^1 \wedge \omega}
\end{align*}
which, restricted to $V$, reads
\begin{equation*}
	0
	= \pa*{\frac{1}{\lambda_2^2} - 2\lambda}\theta_1\, p^1 \wedge \omega + \theta_1\, \rho(a) \pa{p^1 \wedge \omega}.
\end{equation*}
The individual components of $a$ therefore act as
\begin{equation*}
	\rho(a) \pa*{p^1 \wedge \omega}
	= -\frac{3}{2} c\, p^1 \wedge \omega + \pa*{\rho(s) p^1} \wedge \omega + p^1 \wedge \underbrace{\rho(s) \omega}_{= 0}.
\end{equation*}
All together, this leaves us with the scalar equation
\begin{align*}
	0
	&= \theta_1 \cdot \pa*{\frac{1}{\lambda_2^2} - 2\lambda
				- \frac{3}{2} c - \frac{1}{2 \lambda_1^2 \lambda_2^2} A_{11}} \\
	&= \theta_1 \cdot \underbrace{\pa*{\frac{{\beta_{23}}^{2}}{2 \, {\lambda_1} {\lambda_2}}
			+ \frac{{\beta_{34}}^{2}}{2 \, {\lambda_1} {\lambda_2}}
			+ \frac{{\theta_1}^{2}}{2 \, {\lambda_1} {\lambda_2}^{2}}
			+ \frac{1}{2 \, {\lambda_2}^{2}}}}_{> 0},
\end{align*}
contradicting $\theta_1 > 0$.
\end{proof}

The previous lemma shows we can assume $H = \beta \wedge z^*$ for some two-form~$\beta \in \Lambda^2 V^*$.
Together with the simplification~$g = \Diag(\lambda_1, \lambda_2, \lambda_1, \lambda_2, 1)$ from \cref{lem:simplification_metric}, the nilsoliton equations simplify:
\begin{lemma} \label{lem:nil-soliton_red}
On the five-dimensional Heisenberg algebra~$\heis_5$, the nilsoliton equations~\eqref{eq:nil_soliton} with $g = \Diag(\lambda_1, \lambda_2, \lambda_1, \lambda_2, 1)$ and $H = \beta \wedge z^*$ are equivalent to
\begin{equation} \label{eq:nil-soliton_red}
	\left\{
	\begin{aligned}
		0 &= \lambda \tensor{\omega}{_{ij}}
				- 2\tensor{\Ric}{_i^k} \tensor{\omega}{_{kj}}
				+ \frac{1}{2} \tensor{\pa[\big]{\beta^2}}{_{ik}} \tensor{g}{^{kl}} \tensor{\omega}{_{lj}}
				- \frac{1}{2} \tensor{\pa[\big]{\beta^2}}{_{jk}} \tensor{g}{^{kl}} \tensor{\omega}{_{li}}
				+ \frac{1}{4} \pa*{\norm{\omega}_g^2 - \norm{\beta}_g^2} \tensor{\omega}{_{ij}} \\
		0 &= \lambda \tensor{\beta}{_{ij}}
				- \tensor{\Ric}{_i^k} \tensor{\beta}{_{kj}}
				+ \tensor{\Ric}{_j^k} \tensor{\beta}{_{ki}}
				- \frac{1}{4}\pa*{\norm{\omega}_g^2 - \norm{\beta}_g^2} \tensor{\beta}{_{ij}}
				+ \tensor{\pa[\big]{\beta^2}}{_{ik}} \tensor{g}{^{kl}} \tensor{\beta}{_{lj}}
	\end{aligned}
	\right.,
\end{equation}
where
\begin{align}
    \Ric^g_{ij}
    &=
    		\begin{multlined}[t]
    			-\frac{1}{2 \, {\lambda_1}} p^1\otimes p^1
	   			-\frac{1}{2 \, {\lambda_2}} p^2\otimes p^2
	   			-\frac{1}{2 \, {\lambda_1}} q^1\otimes q^1 \\
	   			-\frac{1}{2 \, {\lambda_2}} q^2\otimes q^2
	   			+ \left( \frac{1}{2 \, {\lambda_1}^{2}}
	   			+ \frac{1}{2 \, {\lambda_2}^{2}} \right) z^*\otimes z^*,
    		\end{multlined}
    		\label{eq:ricci_heis2} \\
   	\norm{\omega}_g^2
   	&= \frac{2}{\lambda_1^2} + \frac{2}{\lambda_2^2} \notag
\end{align}
and the symmetric $(0, 2)$-tensor~$(\beta)^2$ for a given two-form~$\beta$ is defined as
\begin{inline}
	\tensor{\pa[\big]{\beta^2}}{_{ij}} \ceq \tensor{g}{^{ab}} \tensor{\beta}{_{ai}} \tensor{\beta}{_{bj}}.
\end{inline}
\end{lemma}
\begin{proof}
The expressions for $\Ric^g$ and $\norm{\omega}_g^2$ follow from a direct calculation.

By the first equation in \labelcref{eq:nil_soliton},
\begin{inline}
	\Ric_g^{1,1} - \frac{1}{4} g^{-1} (\beta \wedge z^*)^2 - \lambda \Id
\end{inline}
is a derivation of $\heis_5$, so
\begin{equation*}
	\rho\pa*{\Ric_g^{1,1} - \frac{1}{4} g^{-1} (\beta \wedge z^*)^2 - \lambda \Id}
			\omega \tprod z
	+ \omega \tprod
			\pa*{\Ric_g^{1,1} - \frac{1}{4} g^{-1} (\beta \wedge z^*)^2 - \lambda \Id}z
	= 0.
\end{equation*}
Simplifying this equation yields
\begin{align*}
	\bra*{\rho\pa*{g^{-1} \pa*{\beta \wedge z^*}^2} \omega}\indices{_{ij}}
	&= -\tensor{\omega}{_{kj}} \tensor{g}{^{kl}} {\pa*{\beta \wedge z^*}^2}\indices{_{li}}
			-\tensor{\omega}{_{ik}} \tensor{g}{^{kl}} {\pa*{\beta \wedge z^*}^2}\indices{_{lj}} \\
	&= -2\tensor{\omega}{_{kj}} \tensor{g}{^{kl}} {\pa*{\beta}^2}\indices{_{li}}
			-2\tensor{\omega}{_{ik}} \tensor{g}{^{kl}} {\pa*{\beta}^2}\indices{_{lj}}.
\end{align*}
The almost complex structure
\begin{equation*}
	J
	\ceq \bra{\tensor{\omega}{_{ij}}} \in \End(V),
	\quad\text{that is,}\quad
	Jp_i = q_i,\, Jq_i = -p_i
\end{equation*}
induces the embedding
\begin{equation*}
	\GL(k, \C) \longhookrightarrow \GL(2k, \R),\quad
	X + i Y \longmapsto
			\begin{bmatrix}
				X & Y \\
				-Y & X
			\end{bmatrix},
\end{equation*}
so block matrices of shape $\Diag(X, X)$ with $X \in \R^{k \times k}$ commute with $J$.
Hence, $\Ric^{1,1}|_V$ and $g|_V$ are complex linear and therefore
\begin{align*}
	\bra*{\rho\pa*{\Ric_g^{1,1}} \omega}\indices{_{ij}}
	= -\tensor{\Ric}{_i^k} \tensor{\omega}{_{kj}}
			-\tensor{\Ric}{_j^k} \tensor{\omega}{_{ik}}
	= -2 \tensor{\Ric}{_i^k} \tensor{\omega}{_{kj}}.
\end{align*}

On the other hand, from
\begin{gather*}
	{\pa*{\beta \wedge z^*}^2}\indices{_{i\mathrm{z}}}
	= \pa*{\beta \wedge z^*}\indices{_{abi}} \pa*{\beta \wedge z^*}\indices{^{ab}_{\mathrm{z}}}
	= \tensor{\beta}{_{ab}} \tensor{\beta}{^{ab}} \tensor{\delta}{_{i\mathrm{z}}},
\shortintertext{we see}
	\frac{1}{4} g^{-1} (\beta \wedge z^*)^2 (z)
	= \frac{1}{4} \norm{\beta}_g^2 \cdot z
\end{gather*}
and from \cref{eq:ricci_heis2}, evidently,
\begin{inline}
	\Ric_g^{1, 1}(z)
	= \frac{1}{4} \norm{\omega}_g^2 \cdot z.
\end{inline}
Together, we obtain the first equation in \labelcref{eq:nil-soliton_red}.
\medbreak

Concerning the second equation of the soliton equations, $\Laplace^g(\beta \wedge z^*) = 0$ yields
\begin{align*}
	0
	&= \lambda \beta \wedge z^*
			+ \rho\pa*{\Ric_g^{1,1} - \frac{1}{4} g^{-1} (\beta \wedge z^*)^2} (\beta \wedge z^*) \\
	&= \lambda \beta \wedge z^*
			+ \rho\pa*{\Ric_g^{1,1} - \frac{1}{4} g^{-1} (\beta \wedge z^*)^2} \beta \wedge z^*
			+ \beta \wedge \rho\pa*{\Ric_g^{1,1} - \frac{1}{4} g^{-1} (\beta \wedge z^*)^2} z^*.
\end{align*}
By the same arguments as above,
\begin{equation*}
	\rho\pa*{\Ric_g^{1,1} - \frac{1}{4} g^{-1} (\beta \wedge z^*)^2} z^*
	= -\frac{1}{4} \pa*{\norm{\omega}_g^2 - \norm{\beta}_g^2} \cdot z^*.
\end{equation*}
The second equation now follows as well, because clearly,
\begin{align*}
	\bra*{\rho\pa*{-\frac{1}{4} g^{-1} (\beta \wedge z^*)^2} \beta}\indices{_{ij}}
	&= \frac{1}{2} (\beta^2)\indices{_{ik}} \tensor{g}{^{kl}} \tensor{\beta}{_{lj}}
			+ \frac{1}{2} (\beta^2)\indices{_{jk}} \tensor{g}{^{kl}} \tensor{\beta}{_{il}} \\
	&= \frac{1}{2} (\beta^2)\indices{_{ik}} \tensor{g}{^{kl}} \tensor{\beta}{_{lj}}
			+ \frac{1}{2} (\tensor{g}{^{ab}}\tensor{\beta}{_{aj}} \tensor{\beta}{_{bk}}) \tensor{g}{^{kl}} \tensor{\beta}{_{il}} \\
	&= (\beta^2)\indices{_{ik}} \tensor{g}{^{kl}} \tensor{\beta}{_{lj}}.
	\qedhere
\end{align*}
\end{proof}

From this, we can fully determine the metric and scaling parameter~$\lambda$ of solitons on $\heis_5$:
\begin{lemma}
For all nilsolitons~$(g, H, \lambda)$ on $\heis_5$ with $g = \Diag(\lambda_1, \lambda_2, \lambda_1, \lambda_2, 1)$ in the ordered basis~$(p_1, p_2, q_1, q_2, z)$, necessarily,
\begin{enumerate*}
	\item $\lambda_1 = \lambda_2$ and
	\item $\lambda = -\frac{2}{\lambda_1^2}$.
\end{enumerate*}
\end{lemma}
\begin{proof}
By \cref{lem:heis5_solitons_harmonic}, $H = \beta \wedge z^*$ for some $\beta \in \Lambda^2 V^*$.
We may take the symplectic trace of the first equation in \labelcref{eq:nil-soliton_red} restricted to the symplectic space~$V$.
This leads to
\begin{align*}
    0
    &= 4\lambda
    			- 2 \Trace(\Ric^{1,1}|_V)
    			+ \norm{\beta}_g^2
    			+ \frac{1}{4} \pa*{\norm{\omega}_g^2 - \norm{\beta}_g^2} \cdot 4 \\
    &= 4\lambda -2 \pa*{-\frac{1}{\lambda_1^2} - \frac{1}{\lambda_2^2}} + 2 \pa*{\frac{1}{\lambda_1^2} + \frac{1}{\lambda_2^2}},
\end{align*}
where the last equality follows from $\norm{\omega}_g^2 = \frac{2}{\lambda_1^2} + \frac{2}{\lambda_2^2}$ and the explicit shape of $\tensor*{\Ric}{^g_{ij}}$.
Thus,
\begin{equation*}
	\lambda
	= -\frac{1}{\lambda_1^2} - \frac{1}{\lambda_2^2}.
\end{equation*}

Concerning the first claim, due to \cref{lem:nil-soliton_red}, the tensor
\begin{equation*}
	\tensor{A}{_{ij}}
	\ceq \lambda \tensor{\omega}{_{ij}}
			- 2\tensor{\Ric}{_i^k} \tensor{\omega}{_{kj}}
			+ \frac{1}{2} \tensor{\pa[\big]{\beta^2}}{_{ik}} \tensor{g}{^{kl}} \tensor{\omega}{_{lj}}
			- \frac{1}{2} \tensor{\pa[\big]{\beta^2}}{_{jk}} \tensor{g}{^{kl}} \tensor{\omega}{_{li}}
			+ \frac{1}{4} \pa*{\norm{\omega}_g^2 - \norm{\beta}_g^2} \tensor{\omega}{_{ij}}
\end{equation*}
has to vanish.
In particular, we calculate the $\tensor{A}{_{13}}$-component and obtain
\begin{equation*}
	0
	= A(p_1, q_1)
	= \frac{\beta(p_1, q_1)^2}{2\lambda_1^2} - \frac{\beta(p_2, q_2)^2}{2\lambda_2^2} + \frac{1}{2\lambda_1^2} - \frac{1}{2\lambda_2^2}.
\end{equation*}
But $H$ is closed, so
\begin{equation*}
    0
    = \D\pa*{\beta \wedge z^*}
    = \pa*{\beta(p_1, q_1) + \beta(p_2, q_2)}\, p^1 \wedge p^2 \wedge q^1 \wedge q^2.
\end{equation*}
All together, this shows $\lambda_1 = \lambda_2$ and $\lambda = -\frac{2}{\lambda_1^2}$.
\end{proof}

Finally, we consider the isotropy group~$\Aut(\heis_5)_g$ of all automorphisms which fix the metric.
Before doing so, we may rescale the metric using \cref{lem:rescaling_soliton} together with the action of $\R^+ \le \Aut(\heis_5)^0$ and without loss of generality assume $[\tensor{g}{_{ij}}] = [\tensor{\delta}{_{ij}}]$ as well as $\lambda = -2$.

\begin{lemma}
The stabilizer~$\Aut(\heis_n)_g$ of all automorphisms that preserve the metric~$g$ given by $[\tensor{g}{_{ij}}] = [\tensor{\delta}{_{ij}}]$ in the canonical basis~\eqref{eq:ordered_basis_heis} with $n = 2k + 1$ is given by
\begin{equation*}
	\UU(k) \rtimes \Z_2
	= \iset{A \in \Aut(\heis_n)}{A^* g = g,\, A^*\omega = \pm \omega},
\end{equation*}
where $\UU(k) \le \Sp(2k, \R) \le \Aut(\heis_n)^0$ and $\Z_2 \le \Aut(\heis_n)$ is generated by $\bar{\Omega}$ from \cref{lem:aut_heisk}.
\end{lemma}
\begin{proof}
It is clear that the $\Z_2$-action preserves~$g$.
By \cref{eq:decomp_auto} and
\begin{equation*}
	\begin{bmatrix}
		S \sqrt{c} & 0 \\
		d^\intercal & c
	\end{bmatrix}^\intercal
			\begin{bmatrix}
				I_{2k} & 0 \\
				0 & 1
			\end{bmatrix}
			\begin{bmatrix}
				S \sqrt{c} & 0 \\
				d^\intercal & c
			\end{bmatrix}
	=
			\begin{bmatrix}
				c \cdot S^\intercal S + d d^\intercal & c d \\
				c d^\intercal & c^2
			\end{bmatrix},
\end{equation*}
in order to preserve~$g$, necessarily $d = 0$ and $c = 1$.
Lastly, $S^\intercal S = I_{2k}$ implies
\begin{equation*}
	S \in \mathrm{O}(2k, \R) \cap \Sp(2k, \R) = \UU(k),
\end{equation*}
which yields the claimed statement.
\end{proof}

This now allows us to prove the second main theorem of this paper:
\begin{theorem} \label{thm:grs_heis5}
Up to automorphism and rescaling, all nilsolitons~$(g, H, \lambda)$ on the 5-dimensional Heisenberg algebra in the ordered basis~$\pa{p_1, p_2, q_1, q_2, z}$ are given by $[g_{ij}] = [\delta_{ij}]$, $\lambda = -2$ and either
\begin{enumerate}
\item $H \ceq 0$ (classical),
\item $H \ceq \pa*{p^1 \wedge q^1 - p^2 \wedge q^2} \wedge z^*$,
\item $H \ceq \pa*{p^1 \wedge p^2 - q^1 \wedge q^2} \wedge z^*$ or
\item $H \ceq \frac{1}{3} \sqrt{3} \cdot \pa*{p^1 \wedge p^2
			+ p^1 \wedge q^1
			- p^2 \wedge q^2
			- q^1 \wedge q^2} \wedge z^*$,
\end{enumerate}
where the Lie bracket is given by $[p_1, q_1] = [p_2, q_2] = z$ and zero otherwise.
None of the above solitons are isomorphic.
\end{theorem}
\begin{proof}
Same as above, we will write $V \ceq z^\perp$ for the orthogonal complement of $z$ with respect to $g$.
By the previous considerations, we can assume $H = \beta \wedge z^*$ for some $\beta \in \Lambda^2 V^*$, $\lambda = -2$ and $[\tensor{g}{_{ij}}] = [\tensor{\delta}{_{ij}}]$.
Recall that $\Aut(\heis_5)_g = \UU(2) \rtimes \Z_2$.
As $\UU(2)$ acts trivially on $z^*$, it is sufficient to study the action of $\UU(2)$ onto $\Lambda^2 V^*$.
In order to do so, we complexify $V$, yielding
\begin{equation} \label{eq:splitting_complex_2}
	\Lambda^2 V^* \le \Lambda_\C^2 V^* \cong \Lambda^{1,1} V^* \oplus \pa*{\Lambda^{0,2} V^* \oplus \Lambda^{2, 0} V^*}.
\end{equation}
It is a well-known fact that $\UU(2) = \pa*{\SU(2) \times \UU(1)}/\set{\pm \Id}$, where $\UU(1) \hookrightarrow \UU(2)$ is embedded via
\begin{inline}
	\alpha \mapsto \Diag(\alpha, \alpha).
\end{inline}
As $\pm\Id$ acts trivially on forms of even degree, it is sensible to study the actions of $\SU(2)$ and $\UU(1)$ on the subspaces given by the splitting in \cref{eq:splitting_complex_2} independently.
Let $\beta = \beta^{1,1} + \beta^{0/2}$ denote the components of $\beta$ with respect to this decomposition.
By definition, $\SU(n)$ fixes holomorphic and anti-holomorphic top-forms, so $\SU(2)$ acts trivially on
\begin{inline}
	\Lambda^{0,2} V^* \oplus \Lambda^{2, 0} V^*.
\end{inline}
Likewise, $U(1)$ acts on $\Lambda^{1,1} V^*$ trivially as well, as for all $\alpha \in \UU(1)$ and $\gamma \in \Lambda^{1,1} V^*$, it is
\begin{inline}
	(\alpha)^* \gamma
	= \alpha \overline{\alpha} \cdot \gamma
	= \gamma.
\end{inline}
\smallbreak

In the canonical complex basis of $V_\C^*$ given by
\begin{align*}
	z^j &\ceq p^j + i p^j \circ J = p^j + i q^j, \\
	\bar{z}^j &\ceq p^j - i p^j \circ J = p^j - i q^j,
\end{align*}
the two-dimensional subspace~$\Lambda^{0,2} V^* \oplus \Lambda^{2, 0} V^* \le \Lambda_\C^2 V^*$ is spanned by the real-valued two-forms
\begin{gather}
	v \ceq z^1 \wedge z^2 + \bar{z}^1 \wedge \bar{z}^2
\quad\text{and}\quad
	w \ceq -i z^1 \wedge z^2 + i \bar{z}^1 \wedge \bar{z}^2.
\end{gather}
Therefore, an arbitrary real two-form in $\Lambda^{0,2} V^* \oplus \Lambda^{2, 0} V^*$ is given by
\begin{inline}
	a v + b w
\end{inline}
with $a, b \in \R$.

It is easy to see that $\UU(1) \cong \SO(2)$ acts via rotation in the (real) plane spanned by $v$ and $w$, (that is, $\R v \oplus \R w$), as
\begin{equation*}
	\pa[\big]{e^{i\varphi}}^*\pa*{av + bw}
	=
			\begin{bmatrix}
				v & w
			\end{bmatrix}
			\begin{bmatrix}
				\cos 2\varphi & \sin 2\varphi \\
				-\sin 2\varphi & \cos 2\varphi
			\end{bmatrix}
			\begin{bmatrix}
				a \\
				b
			\end{bmatrix}.
\end{equation*}
Alas, all non-zero real 2-forms in $\Lambda^{0,2} V^* \oplus \Lambda^{2, 0} V^*$ can be rotated to a real, positive scalar multiple of $v$.
Hence,
\begin{equation*}
	(\alpha)^* \beta^{0/2} = a \cdot v = 2a\, \pa*{p^1 \wedge p^2 - q^1 \wedge q^2}
\end{equation*}
for some $\alpha \in U(1)$ and $a \in \R_0^+$.
\smallbreak

Finally, $\beta^{1,1} \in \Lambda^{1,1} V^* \cap \Lambda^2 V^*$ is skew-hermitian over $(V, J)$, i.e.\ $\beta^{1,1}(J \cdot, J \cdot) = 0$,
so there exists a matrix~$A \in \SU(2)$ that diagonalizes~$\beta^{1,1}$,
\begin{align*}
	A^*\beta^{1,1}
	&= 2c \cdot (Jp^1) \wedge p^1 + 2d \cdot (Jp^2) \wedge p^2 \\
	&= c i \cdot z^1 \wedge \bar{z}^1 + di \cdot z^2 \wedge \bar{z}^2
\end{align*}
for some $c, d \in \R$.
Necessarily, $A^* \beta^{1,1} \wedge z^* \in \Lambda^3 \heis_5^*$ is closed, because $\beta \wedge z^*$ and $\beta^{0/2} \wedge z^*$ are, implying $(A^* \beta^{1,1}) \wedge \omega = 0$.
As
\begin{inline}
	\omega
	= \frac{1}{2}i\, \pa*{z^1 \wedge \bar{z}^1 + z^2 \wedge \bar{z}^2},
\end{inline}
this implies $c = -d$.
A (subsequent) action of
\begin{equation*}
	\begin{bmatrix}
		0 & i \\
		i & 0
	\end{bmatrix}
	\in \SU(2)
\end{equation*}
maps $A^*\beta^{1,1}$ to $-A^*\beta^{1,1}$, so we may assume $c = -d \ge 0$.
The component~$\beta^{1,1}$ is hence put into the shape
\begin{equation*}
	A^*\beta^{1,1} = 2c \pa*{p^1 \wedge q^1 - p^2 \wedge q^2}.
\end{equation*}

\medbreak

In conclusion, we established that every soliton~$(g, H, \lambda)$ on $\heis_5$ up to rescaling is isomorphic to $H = \beta \wedge z^*$ with
\begin{equation*}
	\beta = \tilde{a}\, \pa*{p^1 \wedge p^2 - q^1 \wedge q^2}
			+ \tilde{b}\, \pa*{p^1 \wedge q^1 - p^2 \wedge q^2},\quad \tilde{a}, \tilde{b} \in \R_0^+,
\end{equation*}
together with $[\tensor{g}{_{ij}}] = [\tensor{\delta}{_{ij}}]$ and $\lambda = -2$.

With this, we can now explicitly calculate all quantities in \cref{eq:nil-soliton_red}, yielding
\begin{align*}
	[\Ric_{ij}] &= \Diag\pa*{-\frac{1}{2}, -\frac{1}{2}, -\frac{1}{2}, -\frac{1}{2}, 1}, \\
	\norm{\beta}_g^2 &= 4\tilde{a}^2 + 4\tilde{b}^2, \\
	\norm{\omega}_g^2 &= 4
\shortintertext{and}
	\bra*{\tensor{(\beta)}{^2_{ij}}} &=
			\begin{bmatrix}
				\tilde{a}^{2} + \tilde{b}^{2} & 0 & 0 & 2 \, \tilde{a} \tilde{b} \\
				0 & \tilde{a}^{2} + \tilde{b}^{2} & 2 \, \tilde{a} \tilde{b} & 0 \\
				0 & 2 \, \tilde{a} \tilde{b} & \tilde{a}^{2} + \tilde{b}^{2} & 0 \\
				2 \, \tilde{a} \tilde{b} & 0 & 0 & \tilde{a}^{2} + \tilde{b}^{2}
			\end{bmatrix}.
\end{align*}
Given this, we see that the first equation in \eqref{eq:nil-soliton_red} is now automatically satisfied, while the second one yields the necessary (and sufficient) conditions
\begin{equation*}
	0 = -2 \, {\tilde{a}}^{3} - 4 \, {\tilde{a}} {\tilde{b}}^{2} + 2 \, {\tilde{a}}
\qquad\text{and}\qquad
	0 = -4 \, {\tilde{a}}^{2} {\tilde{b}} - 2 \, {\tilde{b}}^{3} + 2 \, {\tilde{b}}.
\end{equation*}
All solutions are therefore given by either
\begin{enumerate*}
\item $\tilde{a} = \tilde{b} = 0$,
\item $\tilde{a} = 0$ and $\tilde{b} = 1$,
\item $\tilde{a} = 1$ and $\tilde{b} = 0$ or
\item $\tilde{a} = \tilde{b} = \frac{1}{3} \sqrt{3}$.
\end{enumerate*}
\end{proof}

\section{Generalized Algebraic Solitons on the Four-Dimensional Oscillator Algebras} \label{sec:oscillator_algebras}

Next, we investigate algebraic solitons on the four-dimensional spherical and hyperbolic oscillator algebra.
To simplify the treatment of both algebras, we set $\tau \in \set{\pm 1}$.

\begin{definition}
In Salamon notation, the \emph{(spherical) oscillator algebra}~$\mathfrak{osc}_4$ and the \emph{hyperbolic oscillator algebra}~$\mathfrak{osc}'_4$ are the Lie algebras
\begin{equation*}
	\mathfrak{osc}_4 \ceq \mathfrak{osc}_4^{-1} \ceq \pa*{0, +13, -12, -23}
\qquad\text{and}\qquad
	\mathfrak{osc}'_4 \ceq \mathfrak{osc}_4^{1} \ceq \pa*{0, -13, -12, -23}.
\end{equation*}
We denote the ordered base of $\mathfrak{osc}_4$ and $\mathfrak{osc}'_4$ (considered as vector spaces) via $(d, p, q, z)$.
Then equivalently, the Lie bracket on $\mathfrak{osc}_4^\tau$ is given by
\begin{inline}
	[d, p] = q,\, [d, q] = \tau p,\, [p, q] = z
\end{inline}
and zero otherwise.
\end{definition}
It is clear that
\begin{equation*}
	\mathfrak{osc}_4^\tau = \R_d \ltimes_{\operatorname{ad}_d} \heis_3
	\qquad\text{with}\qquad
	\operatorname{ad}_d
	=
			\begin{bmatrix}
				0 & \tau & 0 \\
				1 & 0 & 0 \\
				0 & 0 & 0
			\end{bmatrix}.
\end{equation*}
Also, the ideal~$\heis_3 = [\mathfrak{osc}_4^\tau, \mathfrak{osc}_4^\tau]$ is invariant under automorphisms of the respective oscillator algebras, which allows us to compute the (well-known) automorphism groups
\begin{equation*}
		\Aut(\mathfrak{osc}_4) = \R^3 \rtimes \underbrace{\pa*{\SO(2) \rtimes \R^\times}}_{\mathrm{(I)}}
		\qquad\text{and}\qquad
		\Aut(\mathfrak{osc}'_4) = \R^3 \rtimes \underbrace{\pa*{\pa*{\SO(1, 1) \rtimes \Z_4} \rtimes \R^\times}}_{\mathrm{(II)}},
\end{equation*}
where $\R^3$ acts via
\begin{equation*}
		\R^3 \longto \Aut(\mathfrak{osc}_4^\tau),\quad
				v
				\longmapsto
				\begin{bmatrix}
					1 & 0 & 0 & 0 \\
					-\tau v_1 & 1 & 0 & 0 \\
					v_2 & 0 & 1 & 0 \\
					\frac{1}{2}\pa*{-\tau v_1^2 + v_2^2} + v_3 & v_1 & v_2 & 1
				\end{bmatrix}
\end{equation*}
and the factors~(I) and (II) correspond to pairs~$(\alpha, \beta) \in \set{\pm 1} \times \Aut(\heis_3)$ that intertwine~$\operatorname{ad}_d$, that is,
\begin{equation*}
	\operatorname{ad}_d \circ \beta = \alpha \cdot \beta \circ \operatorname{ad}_d.
\end{equation*}

Analogous to the classification of solitons on the Heisenberg algebra above, we can study the conjugation action of $\Aut(\mathfrak{osc}_4^\tau)$ on positive definite symmetric inner products, i.e.\
\begin{equation*}
	\Aut(\mathfrak{osc}_4^\tau) \curvearrowright \Sym_+^2\pa*{\mathfrak{osc}^\tau_4}^*,\quad (A, g) \longmapsto A^* g.
\end{equation*}
For a fixed metric~$g$ in the ordered basis~$(d, p, q, z)$, we denote the coefficients by
\begin{equation*}
		[g_{ij}]
		= \begin{bmatrix}
			g_{11} & w_1 & w_2 & g_{1n} \\
			w_1 & a & c & v_1 \\
			w_2 & c & b & v_2 \\
			g_{1n} & v_1 & v_2 & g_{nn}
		\end{bmatrix}.
\end{equation*}
We can always assume $g_{1n} = v_1 = v_2 = 0$, because the action of a generic~$\gamma \in \R^3$ on $g$ reads
\begin{multline*}
	\begin{bmatrix}
		1 & 0 & 0 & 0 \\
		-{\gamma_1} {\tau} & 1 & 0 & 0 \\
		{\gamma_2} & 0 & 1 & 0 \\
		{\gamma_3} & {\gamma_1} & {\gamma_2} & 1
	\end{bmatrix}^\intercal
	\begin{bmatrix}
		g_{11} & w_1 & w_2 & g_{1n} \\
		w_1 & a & c & v_1 \\
		w_2 & c & b & v_2 \\
		g_{1n} & v_1 & v_2 & g_{nn}
	\end{bmatrix}
	\begin{bmatrix}
		1 & 0 & 0 & 0 \\
		-{\gamma_1} {\tau} & 1 & 0 & 0 \\
		{\gamma_2} & 0 & 1 & 0 \\
		{\gamma_3} & {\gamma_1} & {\gamma_2} & 1
	\end{bmatrix} \\
	=
		\begin{bsmallmatrix}
			{*} & {*} & {*} & {*} \\
			{*} & {g_{nn}} \gamma_1^{2} + 2 \, v_{1} {\gamma_1} + a & {*} &
					{g_{nn}} {\gamma_1} + v_{1} \\
			{*} & {g_{nn}} {\gamma_1} {\gamma_2} + v_{2} {\gamma_1} + v_{1} {\gamma_2} + c &
					{g_{nn}} \gamma_2^{2} + 2 \, v_{2} {\gamma_2} + b & {g_{nn}} {\gamma_2} + v_{2} \\
			-v_{1} {\gamma_1} {\tau} + v_{2} {\gamma_2} + {g_{nn}} {\gamma_3} + {g_{1n}} &
					{g_{nn}} {\gamma_1} + v_{1} & {g_{nn}} {\gamma_2} + v_{2} & {g_{nn}}
		\end{bsmallmatrix}
\end{multline*}
and the resulting linear equations in $\gamma$,
\begin{equation*}
	\left\{
	\begin{aligned}
	- {g_{1n}} &= -v_{1} {\gamma_1} {\tau} + v_{2} {\gamma_2} + {g_{nn}} {\gamma_3} \\
	-v_{1} &= {g_{nn}} {\gamma_1} \\
	- v_{2} &= {g_{nn}} {\gamma_2}
	\end{aligned}
	\right.
\end{equation*}
always have a solution.

Furthermore, by \cref{lem:rescaling_soliton}, we may rescale~$g$ and assume $g_{11} = 1$.
There also always exists a scaling~$c \in \R^\times \le \Aut(\mathfrak{osc}_4^\tau)$ acting on $g$ such that $g_{nn} = 1$.

Thus, without loss of generality, the metric reads
\begin{equation} \label{eq:partially_simpl_metric}
	[g_{ij}]
	=
	\begin{bmatrix}
		1 & w_1 & w_2 & 0 \\
		w_1 & a & c & 0 \\
		w_2 & c & b & 0 \\
		0 & 0 & 0 & 1
	\end{bmatrix}.
\end{equation}

From this, we can now prove the following theorem:

\begin{theorem} \label{thm:grs_osc}
There do not exist any algebraic solitons (with positive definite background metric~$g$) on the four-dimensional oscillator algebra.
\end{theorem}
\begin{proof}
We can always find a matrix~$S \in O(2)$, such that $S(w_1, w_2)^\intercal = (\sigma, 0)^\intercal$ with $\sigma \ge 0$.

Now, setting up the notation
\begin{align*}
	H
	&= H_{123} d^*\wedge p^*\wedge q^*
			+ H_{124} d^*\wedge p^*\wedge z^* \\
			&\qquad+ H_{134} d^*\wedge q^*\wedge z^*
			+ H_{234} p^*\wedge q^*\wedge z^*
	\in \Lambda^3 \pa{\mathfrak{osc}_4^{-1}}^*,
\end{align*}
we can directly calculate
\begin{align*}
	\det g &= -b {\sigma}^{2} + a b - c^{2}, \\
	[g^{ij}] &=
		\frac{1}{\det g}
		\begin{bmatrix}
			a b - c^{2} & -b {\sigma} & c {\sigma} & 0 \\
			-b {\sigma} & b & -c & 0 \\
			c {\sigma} & -c & -{\sigma}^{2} + a & 0 \\
			0 & 0 & 0 & -b {\sigma}^{2} + a b - c^{2}
		\end{bmatrix}, \\
	[\Ric^g_{ij}] &=
			\frac{1}{2 \det g}
			\begin{bsmallmatrix}
				2 \, {\left(a - b\right)} {\sigma}^{2} - a^{2} + 2 \, a b - b^{2} - 4 \, c^{2} & {\left(a^{2} - b^{2}\right)} {\sigma} & 2 \, {\left(a + b\right)} c {\sigma} & {\sigma}^{2} \\
				{\left(a^{2} - b^{2}\right)} {\sigma} & a^{3} - a b^{2} + 2 \, {\left(a + b\right)} c^{2} + {\sigma}^{2} - a & {\left(a^{2} + 2 \, a b + b^{2} - 1\right)} c & a {\sigma} \\
				2 \, {\left(a + b\right)} c {\sigma} & {\left(a^{2} + 2 \, a b + b^{2} - 1\right)} c & b^{3} + 2 \, {\left(a + b\right)} c^{2} - {\left(a^{2} + 1\right)} b & c {\sigma} \\
				{\sigma}^{2} & a {\sigma} & c {\sigma} & 1
			\end{bsmallmatrix}, \\
	\Laplace H
	&= \left( -\frac{H_{234} {\sigma}^{2} - {\left(H_{134} a - H_{124} c\right)}
					{\sigma} - H_{123}}{b {\sigma}^{2} - a b + c^{2}} \right) d^*\wedge p^*\wedge q^* \\
			&\qquad+ \left( \frac{H_{124} b^{2} + H_{124} c^{2} + H_{234} c {\sigma}
					- {\left(H_{134} a + H_{134} b\right)} c}{b {\sigma}^{2} - a b + c^{2}} \right)
							d^*\wedge p^*\wedge z^* \\
			&\qquad+ \left( \frac{H_{134} a^{2} + H_{134} c^{2}
					- {\left(H_{124} a + H_{124} b\right)} c
					- {\left(H_{234} a - H_{123}\right)}
							{\sigma}}{b {\sigma}^{2} - a b + c^{2}} \right) d^*\wedge q^*\wedge z^*
\end{align*}
and
\begin{align*}
[H^2]_{\bullet1} &= \frac{2}{\det g}
	\begin{bsmallmatrix}
		-H_{134}^{2} {\sigma}^{2} + H_{134}^{2} a + H_{124}^{2} b  - 2 \, H_{124} H_{134} c + H_{123}^{2} \\
		-H_{134} H_{234} {\sigma}^{2} + H_{134} H_{234} a - H_{124} H_{234} c  + {\left(H_{124}^{2} b - H_{124} H_{134} c + H_{123}^{2}\right)} {\sigma} \\
		-H_{124} H_{234} b + H_{134} H_{234} c + {\left(H_{124} H_{134} b - H_{134}^{2} c\right)} {\sigma} \\
		-H_{123} H_{134} {\sigma} + H_{123} H_{234}
	\end{bsmallmatrix}, \\
[H^2]_{\bullet2} &= \frac{2}{\det g}
	\begin{bsmallmatrix}
		-H_{134} H_{234} {\sigma}^{2} + H_{134} H_{234} a - H_{124} H_{234} c  + {\left(H_{124}^{2} b - H_{124} H_{134} c + H_{123}^{2}\right)} {\sigma} \\
		H_{124}^{2} a b - H_{124}^{2} c^{2} - 2 \, H_{124} H_{234} c {\sigma}  - H_{234}^{2} {\sigma}^{2} + {\left(H_{123}^{2} + H_{234}^{2}\right)} a \\
		H_{124} H_{134} a b - H_{124} H_{134} c^{2} + {\left(H_{123}^{2} + H_{234}^{2}\right)} c   - {\left(H_{124} H_{234} b + H_{134} H_{234} c\right)} {\sigma} \\
		-H_{123} H_{134} a + H_{123} H_{124} c + H_{123} H_{234} {\sigma}
	\end{bsmallmatrix}, \\
[H^2]_{\bullet3} &= \frac{2}{\det g}
	\begin{bsmallmatrix}
		-H_{124} H_{234} b + H_{134} H_{234} c + {\left(H_{124} H_{134} b - H_{134}^{2} c\right)} {\sigma} \\
		H_{124} H_{134} a b - H_{124} H_{134} c^{2} + {\left(H_{123}^{2} + H_{234}^{2}\right)} c   - {\left(H_{124} H_{234} b + H_{134} H_{234} c\right)} {\sigma} \\
		-H_{134}^{2} c^{2} - 2 \, H_{134} H_{234} b {\sigma} + {\left(H_{134}^{2} a + H_{123}^{2} + H_{234}^{2}\right)} b \\
		H_{123} H_{124} b - H_{123} H_{134} c
	\end{bsmallmatrix}, \\
[H^2]_{\bullet4} &= \frac{2}{\det g}
	\begin{bsmallmatrix}
		-H_{123} H_{134} {\sigma} + H_{123} H_{234} \\
		-H_{123} H_{134} a + H_{123} H_{124} c + H_{123} H_{234} {\sigma} \\
		H_{123} H_{124} b - H_{123} H_{134} c \\
		H_{134}^{2} a + H_{124}^{2} b - 2 \, H_{124} H_{134} c - 2 \, H_{134} H_{234} {\sigma} + H_{234}^{2}
	\end{bsmallmatrix}.
\end{align*}

From here, we can explicitly calculate the conditions on $H$, $\lambda$ and $g$ imposed by the soliton equations in \cref{eq:nil_soliton}, by setting
\begin{align}
	\label{eq:tensor_FST}
	\mathrm{FST}
	&\ceq \det g \cdot \theta\pa*{g^{-1} \Ric^g - \frac{1}{4} g^{-1} H^2 - \lambda \Id} \mu
\shortintertext{and}
	\label{eq:tensor_SND}
	\mathrm{SND}
	&\ceq \det g \cdot \pa*{\Laplace^g H - \lambda H - \rho\pa*{g^{-1} \Ric^g - \frac{1}{4} g^{-1} H^2} H}
\end{align}
to zero.
The components of the first tensor thus translate to the conditions
\begin{enumerate}[label=(\alph*)]
\item \label{it:cond1} $0 = \frac{1}{2} \, H_{123} H_{234}$,%
\item $0 = \frac{1}{2} \, H_{123} H_{124}$,%
\item \label{it:cond2} $0 = \frac{1}{2} \, H_{123} H_{134} + \frac{1}{2} \, {\sigma}$,%
\item $0 = \frac{1}{2} \, H_{124} H_{234} b - \frac{1}{2} \, H_{134} H_{234} c$,%
\item \label{it:cond3} $0 = \frac{1}{2} \, H_{134} H_{234} a - \frac{1}{2} \, H_{124} H_{234} c - \frac{1}{2} \, H_{234}^{2} {\sigma} - \frac{1}{2} \, {\sigma}$,%
\item $0 = -\frac{1}{2} \, H_{124} H_{134} a - \frac{1}{2} \, H_{124} H_{134} b + \frac{1}{2} \, H_{124}^{2} c + \frac{1}{2} \, H_{134}^{2} c + \frac{1}{2} \, H_{124} H_{234} {\sigma} + 2 \, a c + 2 \, b c$,%
\item $0 = -\frac{1}{2} \, H_{134} H_{234} a + \frac{1}{2} \, H_{124} H_{234} c + \frac{1}{2} \, H_{234}^{2} {\sigma} + \frac{1}{2} \, H_{123} H_{134} + {\sigma}$,%
\item $0 = -\frac{1}{2} \, H_{123} H_{124} b + \frac{1}{2} \, H_{123} H_{134} c - \frac{1}{2} \, H_{124} H_{134} {\sigma} + \frac{1}{2} \, H_{124} H_{234} + \frac{5}{2} \, c {\sigma}$,%
\item $0 = -b {\lambda} {\sigma}^{2} + H_{124}^{2} b - H_{124} H_{134} c + a b {\lambda} - c^{2} {\lambda} + \frac{1}{2} \, H_{123}^{2} - \frac{1}{2} \, a^{2} - a b + \frac{3}{2} \, b^{2} + 2 \, c^{2}$,%
\item $0 = -\frac{1}{2} \, H_{123} H_{134} a + \frac{1}{2} \, H_{123} H_{124} c - \frac{1}{2} \, H_{134}^{2} {\sigma} + \frac{1}{2} \, H_{123} H_{234} {\sigma} + \frac{1}{2} \, H_{134} H_{234} - \frac{3}{2} \, a {\sigma} + b {\sigma}$,%
\item $0 = -b {\lambda} {\sigma}^{2} + a b {\lambda} - c^{2} {\lambda} - \frac{1}{2} \, H_{134} H_{234} {\sigma} + H_{123}^{2} + \frac{1}{2} \, H_{234}^{2} + \frac{3}{2}$ and%
\item $0 = b {\lambda} {\sigma}^{2} - H_{134}^{2} a + H_{124} H_{134} c - a b {\lambda} + c^{2} {\lambda} + H_{134} H_{234} {\sigma} - \frac{1}{2} \, H_{123}^{2} - \frac{3}{2} \, a^{2} + a b + \frac{1}{2} \, b^{2} - 2 \, c^{2}$.%
\end{enumerate}

If it was $\sigma \ne 0$, by condition~\labelcref{it:cond3}, $H_{234} \ne 0$.
Therefore, from condition~\labelcref{it:cond1}, $H_{123} = 0$, which contradicts condition~\labelcref{it:cond2}.
Thus, necessarily, $\sigma = 0$.

For this reason, by a subsequent rotation~$S \in O(2)$ we can diagonalize the metric, i.e.\ assume $c = 0$.
All together, the non-redundant equations (in no particular order) now read
\begin{enumerate}[label=(\alph*')]
\item \label{it:cond_red1} $0 = \frac{1}{2} \, H_{123} H_{124}$, %
\item \label{it:cond_red2} $0 = \frac{1}{2} \, H_{123} H_{134}$, %
\item \label{it:cond_red3} $0 = \frac{1}{2} \, H_{123} H_{234}$, %
\item \label{it:cond_red4} $0 = \frac{1}{2} \, H_{134} H_{234} a$, %
\item \label{it:cond_red5} $0 = \frac{1}{2} \, H_{124} H_{234} b$, %
\item \label{it:cond_red6} $0 = -\frac{1}{2} \, H_{124} H_{134} a - \frac{1}{2} \, H_{124} H_{134} b = -\frac{1}{2} \, H_{124} H_{134} (a + b)$, %
\item $0 = H_{124}^{2} b + a b {\lambda} + \frac{1}{2} \, H_{123}^{2} - \frac{1}{2} \, a^{2} - a b + \frac{3}{2} \, b^{2}$, %
\item $0 = -H_{134}^{2} a - a b {\lambda} - \frac{1}{2} \, H_{123}^{2} - \frac{3}{2} \, a^{2} + a b + \frac{1}{2} \, b^{2}$, %
\item $0 = a b {\lambda} + H_{123}^{2} + \frac{1}{2} \, H_{234}^{2} + \frac{3}{2}$, %
\end{enumerate}
for the first tensor and
\begin{enumerate}[label=(\alph*'), resume]
\item $0 = -\frac{3}{2} \, H_{123} H_{134}^{2} a - \frac{3}{2} \, H_{123} H_{124}^{2} b - H_{123} a b {\lambda} - \frac{3}{2} \, H_{123}^{3} - \frac{3}{2} \, H_{123} H_{234}^{2} - \frac{1}{2} \, H_{123} a^{2} + H_{123} a b - \frac{1}{2} \, H_{123} b^{2} - 2 \, H_{123}$, %
\item $0 = -\frac{3}{2} \, H_{124} H_{134}^{2} a - \frac{3}{2} \, H_{124}^{3} b - H_{124} a b {\lambda} - \frac{3}{2} \, H_{123}^{2} H_{124} - \frac{3}{2} \, H_{124} H_{234}^{2} + H_{124} a b - 2 \, H_{124} b^{2}$, %
\item $0 = -\frac{3}{2} \, H_{134}^{3} a - \frac{3}{2} \, H_{124}^{2} H_{134} b - H_{134} a b {\lambda} - \frac{3}{2} \, H_{123}^{2} H_{134} - \frac{3}{2} \, H_{134} H_{234}^{2} - 2 \, H_{134} a^{2} + H_{134} a b$, %
\item $0 = -\frac{3}{2} \, H_{134}^{2} H_{234} a - \frac{3}{2} \, H_{124}^{2} H_{234} b - H_{234} a b {\lambda} - \frac{3}{2} \, H_{123}^{2} H_{234} - \frac{3}{2} \, H_{234}^{3} - \frac{1}{2} \, H_{234}$, %
\end{enumerate}
for the second one.

Now, for \labelcref{it:cond_red1,it:cond_red2,it:cond_red3,it:cond_red4,it:cond_red5,it:cond_red6} to hold, one of the cases
\begin{itemize}
\item $H_{123} = 0, H_{124} = 0, H_{134} = 0$, %
\item $H_{123} = 0, H_{124} = 0, H_{234} = 0$, %
\item $H_{123} = 0, H_{134} = 0, H_{234} = 0$ or %
\item $H_{124} = 0, H_{134} = 0, H_{234} = 0$ %
\end{itemize}
has be satisfied.
We will now go through all these possible cases and show that in none of those, a solution exists.
To simplify the calculations, we will substitute $a + b \to 2 t_1$ and $a - b \to 2 t_2$ and then show all these cases lead to a contradiction.

\paragraph*{Case $H_{123} = H_{124} = H_{134} = 0$:}

In this case, the equations read
\begin{enumerate}
\item $0 = -t_{1}^{2} {\lambda} + t_{2}^{2} {\lambda} - 4 \, t_{1} t_{2} - 2 \, t_{2}^{2}$, %
\item $0 = t_{1}^{2} {\lambda} - t_{2}^{2} {\lambda} - 4 \, t_{1} t_{2} + 2 \, t_{2}^{2}$, %
\item $0 = -H_{234} t_{1}^{2} {\lambda} + H_{234} t_{2}^{2} {\lambda} - \frac{3}{2} \, H_{234}^{3} - \frac{1}{2} \, H_{234}$, %
\item $0 = t_{1}^{2} {\lambda} - t_{2}^{2} {\lambda} + \frac{1}{2} \, H_{234}^{2} + \frac{3}{2}$. %
\end{enumerate}
The first together with the second equation shows $t_1 t_2 = 0$ and $\pa*{t_1^2 - t_2^2} \lambda = -2 t_2^2$.
Plugging this into the fourth equation, we see
\begin{inline}
	2 t_2^2 = \frac{1}{2} \, H_{234}^{2} + \frac{3}{2} > 0,
\end{inline}
so $t_2 \ne 0$, which yields $t_1 = 0$.
However, $t_1 \le 0$ contradicts $g$ being positive definite.

\paragraph*{Case $H_{123} = H_{124} = H_{234} = 0$:}

In this case, the equations read
\begin{enumerate}
\item $0 = t_{1}^{2} {\lambda} - t_{2}^{2} {\lambda} - 4 \, t_{1} t_{2} + 2 \, t_{2}^{2}$, %
\item $0 = t_{1}^{2} {\lambda} - t_{2}^{2} {\lambda} + \frac{3}{2}$, %
\item $0 = -H_{134}^{2} t_{1} - H_{134}^{2} t_{2} - t_{1}^{2} {\lambda} + t_{2}^{2} {\lambda} - 4 \, t_{1} t_{2} - 2 \, t_{2}^{2}$, %
\item $0 = -\frac{3}{2} \, H_{134}^{3} t_{1} - \frac{3}{2} \, H_{134}^{3} t_{2} - H_{134} t_{1}^{2} {\lambda}
		+ H_{134} t_{2}^{2} {\lambda} - H_{134} t_{1}^{2} - 4 \, H_{134} t_{1} t_{2} - 3 \, H_{134} t_{2}^{2}$. %
\end{enumerate}

If $H_{134} = 0$, then (i)+(iii) yields $t_1 t_2 = 0$.
But this the this means that (i)-(ii) reads $0 = -\frac{3}{2} + 2 t_2^2$, so $t_2 \ne 0$.
Therefore, $t_1 = 0$.

If $H_{134} \ne 0$, then $\mathrm{(iii)} - H_{134}^{-1} \cdot \mathrm{(iv)}$ shows
\begin{inline}
	0
	= \frac{1}{2} H_{134}^2 (t_1 + t_2) + t_1^2 + t_2^2,
\end{inline}
which cannot happen, as $t_1 + t_2 = a > 0$.

\paragraph*{Case $H_{123} = H_{134} = H_{234} = 0$:}

In this case, the equations read
\begin{enumerate}
\item $0 = -t_{1}^{2} {\lambda} + t_{2}^{2} {\lambda} - 4 \, t_{1} t_{2} - 2 \, t_{2}^{2}$, %
\item $0 = t_{1}^{2} {\lambda} - t_{2}^{2} {\lambda} + \frac{3}{2}$, %
\item $0 = -\frac{3}{2} \, H_{124}^{3} t_{1} + \frac{3}{2} \, H_{124}^{3} t_{2} - H_{124} t_{1}^{2} {\lambda} + H_{124} t_{2}^{2} {\lambda} - H_{124} t_{1}^{2} + 4 \, H_{124} t_{1} t_{2} - 3 \, H_{124} t_{2}^{2}$, %
\item $0 = H_{124}^{2} t_{1} - H_{124}^{2} t_{2} + t_{1}^{2} {\lambda} - t_{2}^{2} {\lambda} - 4 \, t_{1} t_{2} + 2 \, t_{2}^{2}$. %
\end{enumerate}

Analogous to the previous case,
if $H_{124} = 0$, then (i)+(iv) yields $t_1 t_2 = 0$.
This means that (i)+(ii) reads $0 = \frac{3}{2} - 2 t_2^2$, so in particular $t_2 \ne 0$.
Therefore, $t_1 = 0$.

If $H_{124} \ne 0$, then $-\mathrm{(iv)} - H_{124}^{-1} \cdot \mathrm{(iii)}$ shows
\begin{inline}
	0
	= \frac{1}{2} H_{124}^2 (t_1 - t_2) + t_1^2 + t_2^2,
\end{inline}
which cannot happen, as $t_1 - t_2 = b > 0$.

\paragraph*{Case $H_{124} = H_{134} = H_{234} = 0$:}

In this case, the equations read
\begin{enumerate}
\item $0 = t_{1}^{2} {\lambda} - t_{2}^{2} {\lambda} + \frac{1}{2} \, H_{123}^{2} - 4 \, t_{1} t_{2} + 2 \, t_{2}^{2}$, %
\item $0 = -H_{123} t_{1}^{2} {\lambda} + H_{123} t_{2}^{2} {\lambda} - \frac{3}{2} \, H_{123}^{3} - 2 \, H_{123} t_{2}^{2} - 2 \, H_{123}$, %
\item $0 = t_{1}^{2} {\lambda} - t_{2}^{2} {\lambda} + H_{123}^{2} + \frac{3}{2}$, %
\item $0 = -t_{1}^{2} {\lambda} + t_{2}^{2} {\lambda} - \frac{1}{2} \, H_{123}^{2} - 4 \, t_{1} t_{2} - 2 \, t_{2}^{2}$. %
\end{enumerate}

Now, (i)+(iv) yield $t_1 t_2 = 0$.
But (iii)-(i) shows
\begin{inline}
	0 = \frac{1}{2} H_{123}^2 + \frac{3}{2} - 2\, t_2^2
\end{inline}
so $t_2 \ne 0$, which means $t_1 = 0$.
This again contradicts $g$ being positive definite.
\smallbreak

In conclusion, there exist no solutions to nilsoliton equations.
\end{proof}

The same approach as in the previous theorem also works for the hyperbolic oscillator algebra.
We obtain the following classification of algebraic solitons.
Note, that this implies that the torsion of all solitons on $\mathfrak{osc}_4'$ is harmonic.

\begin{theorem} \label{thm:grs_split_osc}
Up to automorphism and rescaling of the metric, all algebraic solitons~$(g, H, \lambda)$ on the four-dimensional hyperbolic oscillator algebra in the ordered basis~$\pa{d, p, q, z}$ are given by $\lambda = -2$
and either
\begin{enumerate}
\item
\begin{inline}
	[g_{ij}]
		=
			\Diag\pa*{1,\, \frac{1}{2} \sqrt{3},\, \frac{1}{2} \sqrt{3},\, 1},
\end{inline}
$H = 0$ (classical),
\item $[g_{ij}] = [\delta_{ij}]$, $H = +p^* \wedge q^* \wedge z^*$ or
\item $[g_{ij}] = [\delta_{ij}]$, $H = -p^* \wedge q^* \wedge z^*$,
\end{enumerate}
where the Lie bracket is given by
\begin{inline}
	[d, p] = q,\, [d, q] = p,\, [p, q] = z
\end{inline}
and zero otherwise.
None of these three solitons are isomorphic.
\end{theorem}
\begin{proof}
Recall the simplified version of the metric in \cref{eq:partially_simpl_metric}.
Because rotations by automorphisms of $\mathfrak{osc}_4'$ in the $p$-$q$-plane are Lorentzian,
there are now three distinct cases to consider, depending on $w = (w_1, w_2)^\intercal \in \R^2$:
Either $w$ is zero, light-like or space/time-like with respect to the canonical Minkowski metric~$\D q^2 - \D p^2$ on $\operatorname{Span}\set{p, q}$.

\paragraph*{Light-like case:}

In this case, the vector~$w$ takes the shape~$w = (\sigma, \pm \sigma)$ with $\sigma \in \R^\times$.
A generic matrix in the special Lorentzian group~$\SO(1, 1)$ acts on light-like vectors via scaling, i.e.
\begin{equation*}
	\begin{bmatrix}
		a & b \\
		b & a
	\end{bmatrix}
	\begin{bmatrix}
		\sigma \\ \pm\sigma
	\end{bmatrix}
	= (a \pm b)
			\begin{bmatrix}
				\sigma \\ \pm\sigma
			\end{bmatrix}.
\end{equation*}
Consequently, a hyperbolic rotation that rescales $w$ to $\tilde{w} = (1, \pm 1)^\intercal$ is given by $a \ceq \frac{1}{2}\pa*{\sigma^{-1} + \sigma}$ and $b \ceq \pm\frac{1}{2}\pa*{\sigma^{-1} - \sigma}$.
By a subsequent action of $\Z_4 \le \Aut(\mathfrak{osc}_4')$ we can further assume $\tilde{w} = (1, 1)^\intercal$.

The generic reduced metric now takes the shape\footnote{We chose the center block in this particular way to simplify the calculations later (specifically, this minimizes the number of symbolic summands of the determinant).}
\begin{equation*}
	[g_{ij}]
	=
			\begin{bmatrix}
				1 & 1 & 1 & 0 \\
				1 & 1 + a & 1 - c & 0 \\
				1 & 1 - c & 1 + b & 0 \\
				0 & 0 & 0 & 1
			\end{bmatrix}.
\end{equation*}
Accordingly,
\begin{align*}
	\det g &= a b - c^{2}, \\
	[g^{ij}] &=
		\frac{1}{\det g}
		\begin{bmatrix}
			\det g + a + b + 2 \, c & -b - c & -a - c & 0 \\
			-b - c & b & c & 0 \\
			-a - c & c & a & 0 \\
			0 & 0 & 0 & \det g
		\end{bmatrix}, \\
	\Ric^g_{\bullet1} &= \frac{1}{2 \det g}
		\begin{bsmallmatrix}
			-(a + b)^2 + 4 \, c^{2} \\
			a^{2} - b^{2} + 2 \, {\left(a - b\right)} c \\
			-a^{2} + b^{2} - 2 \, {\left(a - b\right)} c \\
			0
		\end{bsmallmatrix}, \\
	\Ric^g_{\bullet2} &= \frac{1}{2 \det g}
		\begin{bsmallmatrix}
			a^{2} - b^{2} + 2 \, {\left(a - b\right)} c \\
			a^{3} - {\left(a + 1\right)} b^{2}  - 2 \, {\left(a - b\right)} c^{2} + 3 \, a^{2}  - 2 \, a b + 4 \, {\left(a - b\right)} c - a \\
			(a - b)^2 - (a - b)^2 c + c \\
			-a - c
		\end{bsmallmatrix}, \\
	\Ric^g_{\bullet3} &= \frac{1}{2 \det g}
		\begin{bsmallmatrix}
			-a^{2} + b^{2} - 2 \, {\left(a - b\right)} c \\
			(a - b)^2 - (a - b)^2 c + c \\
			b^{3} + 2 \, {\left(a - b\right)} c^{2} - a^{2}  - (a + 1)^2 b + 3 \, b^{2}  - 4 \, {\left(a - b\right)} c \\
			b + c
		\end{bsmallmatrix}, \\
	\Ric^g_{\bullet4} &= \frac{1}{2 \det g}
		\begin{bsmallmatrix}
			0 \\
			-a - c \\
			b + c \\
			1
		\end{bsmallmatrix}
\end{align*}
and
\begin{align*}
[H^2]_{\bullet1} &= \frac{2}{\det g}
	\begin{bsmallmatrix}
		{H_{134}}^{2} a + {H_{124}}^{2} b + 2 \, {H_{124}} {H_{134}} c + {H_{123}}^{2} \\
		{H_{124}}^{2} b + {H_{123}}^{2} + {\left({H_{124}} {H_{134}} + {H_{134}} {H_{234}}\right)} a + {\left({H_{124}}^{2} + {H_{124}} {H_{134}} + {H_{124}} {H_{234}}\right)} c \\
		{H_{134}}^{2} a + {H_{123}}^{2} + {\left({H_{124}} {H_{134}} - {H_{124}} {H_{234}}\right)} b + {\left({H_{124}} {H_{134}} + {H_{134}}^{2} - {H_{134}} {H_{234}}\right)} c \\
		{H_{123}} {\left({H_{124}} - {H_{134}} + {H_{234}}\right)}
	\end{bsmallmatrix}, \\
[H^2]_{\bullet2} &= \frac{2}{\det g}
	\begin{bsmallmatrix}
		{H_{124}}^{2} b + {H_{123}}^{2} + {\left({H_{124}} {H_{134}} + {H_{134}} {H_{234}}\right)} a + {\left({H_{124}}^{2} + {H_{124}} {H_{134}} + {H_{124}} {H_{234}}\right)} c \\
		\mathrm{(III)} \\
		\mathrm{(IV)} \\
		-{\left({H_{134}} a + {H_{124}} c - {H_{124}} + {H_{134}} - {H_{234}}\right)} {H_{123}}
	\end{bsmallmatrix}, \\
[H^2]_{\bullet3} &= \frac{2}{\det g}
	\begin{bsmallmatrix}
		{H_{134}}^{2} a + {H_{123}}^{2} + {\left({H_{124}} {H_{134}} - {H_{124}} {H_{234}}\right)} b + {\left({H_{124}} {H_{134}} + {H_{134}}^{2} - {H_{134}} {H_{234}}\right)} c \\
		\mathrm{(IV)} \\
		\mathrm{(V)} \\
		{\left({H_{124}} b + {H_{134}} c + {H_{124}} - {H_{134}} + {H_{234}}\right)} {H_{123}}
	\end{bsmallmatrix}, \\
[H^2]_{\bullet4} &= \frac{2}{\det g}
	\begin{bsmallmatrix}
		{H_{123}} {\left({H_{124}} - {H_{134}} + {H_{234}}\right)} \\
		-{\left({H_{134}} a + {H_{124}} c - {H_{124}} + {H_{134}} - {H_{234}}\right)} {H_{123}} \\
		{\left({H_{124}} b + {H_{134}} c + {H_{124}} - {H_{134}} + {H_{234}}\right)} {H_{123}} \\
		{H_{134}}^{2} a + {H_{124}}^{2} b + 2 \, {H_{124}} {H_{134}} c + {H_{124}}^{2} - 2 \, {H_{124}} {H_{134}} + {H_{134}}^{2} + 2 \, {\left({H_{124}} - {H_{134}}\right)} {H_{234}} + {H_{234}}^{2}
	\end{bsmallmatrix},
\end{align*}
where
\begin{align*}
	\mathrm{(III)} &\ceq
			-{H_{124}}^{2} c^{2} + {H_{123}}^{2} + {\left({H_{123}}^{2} + {H_{124}}^{2} + 2 \, {H_{124}} {H_{234}} + {H_{234}}^{2}\right)} a \\
			&\qquad+ {\left({H_{124}}^{2} a + {H_{124}}^{2}\right)} b + 2 \, {\left({H_{124}}^{2} + {H_{124}} {H_{234}}\right)} c, \\
	\mathrm{(IV)} &\ceq
			-{H_{124}} {H_{134}} c^{2}
			+ {H_{123}}^{2}
			+ {H_{134}} {\left({H_{124}}  + {H_{234}}\right)} a
			+ {H_{124}} {\left({H_{134}} a + {H_{134}} - {H_{234}}\right)} b \\
			&\qquad- {\left({H_{123}}^{2} - 2 \, {H_{124}} {H_{134}}
			+ {\left({H_{124}} - {H_{134}}\right)} {H_{234}} + {H_{234}}^{2}\right)} c, \\
	\mathrm{(V)} &\ceq
			-{H_{134}}^{2} c^{2} + {H_{134}}^{2} a + {H_{123}}^{2} + 2 \, {\left({H_{134}}^{2} - {H_{134}} {H_{234}}\right)} c \\
			&\qquad+ {\left({H_{134}}^{2} a + {H_{123}}^{2} + {H_{134}}^{2} - 2 \, {H_{134}} {H_{234}} + {H_{234}}^{2}\right)} b.
\end{align*}
From here, we can calculate the tensor~$\operatorname{FST}$, as defined in \cref{eq:tensor_FST}.
By setting its components to zero, we obtain necessary conditions for $(g, H, \lambda)$ to define an algebraic soliton.
Among those,
\begin{enumerate}
\item $0 = -\frac{1}{2} \, {H_{123}} {H_{234}}$, %
\item $0 = \frac{1}{2} \, {H_{123}} {H_{124}} - \frac{1}{2}$ and %
\item $0 = \frac{1}{2} \, {H_{124}} {H_{234}} b + \frac{1}{2} \, {H_{134}} {H_{234}} c + \frac{1}{2} \, {H_{124}} {H_{234}} - \frac{1}{2} \, {H_{134}} {H_{234}} + \frac{1}{2} \, {H_{234}}^{2} + \frac{1}{2}$. %
\end{enumerate}
But these cannot be satisfied simultaneously for the following reason:
The third equation implies $H_{234} \ne 0$ so by the first equation, $H_{123} = 0$, which, however, contradicts the second equation.
Alas, there exist no solutions in the light-like case.

We will now, by a very similar argument, show that in the space-like and time-like case no solutions exist:

\paragraph{Space and time-like case:}

In this case, there exists a hyperbolic rotation~$S \in \SO(1, 1) \rtimes \Z_4$ together with a flip, to such that $\tilde{w} = Sw = (\sigma, 0)^\intercal$, where $\sigma \in \R^+$.

Similar to the light-like case, the generic reduced metric now takes the shape
\begin{equation*}
	[g_{ij}]
	=
			\begin{bmatrix}
				1 & \sigma & 0 & 0 \\
				\sigma & a & c & 0 \\
				0 & c & b & 0 \\
				0 & 0 & 0 & 1
			\end{bmatrix}.
\end{equation*}
Consequently,
\begin{align*}
	\det g &= -b \sigma^2 + a b - c^{2}, \\
	[g^{ij}] &=
		\frac{1}{\det g}
		\begin{bmatrix}
			a b - c^{2} & -b {\sigma} & c {\sigma} & 0 \\
			-b {\sigma} & b & -c & 0 \\
			c {\sigma} & -c & -{\sigma}^{2} + a & 0 \\
			0 & 0 & 0 & \det g
		\end{bmatrix}, \\
	[\Ric^g_{ij}] &= \frac{1}{2 \det g}
	        \begin{bsmallmatrix}
				2 \, {\left(a + b\right)} {\sigma}^{2}  - (a + b)^2 + 4 \, c^{2} & {\left(a^{2} - b^{2}\right)} {\sigma} & 2 \, {\left(a - b\right)} c {\sigma} & -{\sigma}^{2} \\
				{\left(a^{2} - b^{2}\right)} {\sigma} & a^{3} - a b^{2} - 2 \, {\left(a - b\right)} c^{2}  + {\sigma}^{2} - a & (a - b)^2 c - c & -a {\sigma} \\
				2 \, {\left(a - b\right)} c {\sigma} & (a - b)^2 c - c & b^{3} + 2 \, {\left(a - b\right)} c^{2}  - {\left(a^{2} + 1\right)} b & -c {\sigma} \\
				-{\sigma}^{2} & -a {\sigma} & -c {\sigma} & 1
	        \end{bsmallmatrix}
\end{align*}
and the columns of $H^2$ are given by are given by
\begin{align*}
[H^2]_{\bullet1} &= \frac{2}{\det g}
	\begin{bmatrix}
		-H_{134}^{2} {\sigma}^{2} + H_{134}^{2} a + H_{124}^{2} b - 2 \, H_{124} H_{134} c + H_{123}^{2} \\
		-H_{134} H_{234} {\sigma}^{2} + H_{134} H_{234} a - H_{124} H_{234} c + {\left(H_{124}^{2} b - H_{124} H_{134} c + H_{123}^{2}\right)} {\sigma} \\
		-H_{124} H_{234} b + H_{134} H_{234} c + {\left(H_{124} H_{134} b - H_{134}^{2} c\right)} {\sigma} \\
		-H_{123} H_{134} {\sigma} + H_{123} H_{234}
	\end{bmatrix}, \\
[H^2]_{\bullet2} &= \frac{2}{\det g}
	\begin{bmatrix}
		-H_{134} H_{234} {\sigma}^{2} + H_{134} H_{234} a - H_{124} H_{234} c + {\left(H_{124}^{2} b - H_{124} H_{134} c + H_{123}^{2}\right)} {\sigma} \\
		-2 \, H_{124} H_{234} c {\sigma} - H_{234}^{2} {\sigma}^{2} + H_{124}^{2} {{\det g}} + {\left(H_{123}^{2} + H_{234}^{2}\right)} a \\
		H_{124} H_{134} {{\det g}} + {\left(H_{123}^{2} + H_{234}^{2}\right)} c - {\left(H_{124} H_{234} b + H_{134} H_{234} c\right)} {\sigma} \\
		-H_{123} H_{134} a + H_{123} H_{124} c + H_{123} H_{234} {\sigma}
	\end{bmatrix}, \\
[H^2]_{\bullet3} &= \frac{2}{\det g}
	\begin{bmatrix}
		-H_{124} H_{234} b + H_{134} H_{234} c + {\left(H_{124} H_{134} b - H_{134}^{2} c\right)} {\sigma} \\
		H_{124} H_{134} {{\det g}} + {\left(H_{123}^{2} + H_{234}^{2}\right)} c - {\left(H_{124} H_{234} b + H_{134} H_{234} c\right)} {\sigma} \\
		-2 \, H_{134} H_{234} b {\sigma} + H_{134}^{2} {{\det g}} + {\left(H_{123}^{2} + H_{234}^{2}\right)} b \\
		H_{123} H_{124} b - H_{123} H_{134} c
	\end{bmatrix}, \\
[H^2]_{\bullet4} &= \frac{2}{\det g}
	\begin{bmatrix}
		-H_{123} H_{134} {\sigma} + H_{123} H_{234} \\
		-H_{123} H_{134} a + H_{123} H_{124} c + H_{123} H_{234} {\sigma} \\
		H_{123} H_{124} b - H_{123} H_{134} c \\
		H_{134}^{2} a + H_{124}^{2} b - 2 \, H_{124} H_{134} c - 2 \, H_{134} H_{234} {\sigma} + H_{234}^{2}
	\end{bmatrix}.
\end{align*}

Completely analogously to the light-like case, from there, we may calculate the $\operatorname{FST}$-tensor as defined in \cref{eq:tensor_FST} and get a contradiction from the conditions
\begin{enumerate}
\item $0 = -\frac{1}{2} \, H_{123} H_{234}$, %
\item $0 = -\frac{1}{2} \, H_{134} H_{234} a + \frac{1}{2} \, H_{124} H_{234} c + \frac{1}{2} \, H_{234}^{2} {\sigma} + \frac{1}{2} \, {\sigma}$ and %
\item $0 = -\frac{1}{2} \, H_{123} H_{134} + \frac{1}{2}\, \sigma$. %
\end{enumerate}

In conclusion, in order for $g$ to solve the soliton equations, necessarily, $w = 0$ in \cref{eq:partially_simpl_metric}.

\paragraph{Zero-case:}

If $w = 0$, according to \cref{prop:diagonalization_hyperbolic_2x2}, we can choose a hyperbolic rotation to fully diagonalize the metric.
It therefore takes the shape $[g_{ij}] = \Diag(1, a, b, 1)$ for $a, b \in \R^+$.
Analogously to the cases above, we calculate
\begin{align*}
	\det g &= ab, \\
	[g^{ij}] &= \Diag\pa*{1, a^{-1}, b^{-1}, 1}, \\
	[\Ric^g_{ij}] &= \frac{1}{2ab}
			\Diag\pa*{
				-(a + b)^2,\,
				(a^{2} - b^{2} - 1)a,\,
				-(a^{2} - b^{2} + 1)b,\,
			1}, \\
	\Laplace^g H &= -\frac{H_{123}}{a b} d^*\wedge p^*\wedge q^* -\frac{H_{124} b}{a} d^*\wedge p^*\wedge z^* -\frac{H_{134} a}{b} d^*\wedge q^*\wedge z^*, \\
	H^2 &= [H^2_{1,2}, H^2_{3,4}]
\end{align*}
where $H^2$ is comprised of the block matrices
\begin{align*}
	H^2_{1,2} &= \frac{2}{ab}
			\begin{bsmallmatrix}
				H_{134}^{2} a + H_{124}^{2} b + H_{123}^{2} & H_{134} H_{234} a \\
				H_{134} H_{234} a & {\left(H_{124}^{2} b + H_{123}^{2} + H_{234}^{2}\right)} a \\
				-H_{124} H_{234} b & H_{124} H_{134} a b \\
				H_{123} H_{234} & -H_{123} H_{134} a
			\end{bsmallmatrix}, \\
	H^2_{3,4} &= \frac{2}{ab}
			\begin{bsmallmatrix}
				-H_{124} H_{234} b & H_{123} H_{234} \\
				H_{124} H_{134} a b & -H_{123} H_{134} a \\
				{\left(H_{134}^{2} a + H_{123}^{2} + H_{234}^{2}\right)} b & H_{123} H_{124} b \\
				H_{123} H_{124} b & H_{134}^{2} a + H_{124}^{2} b + H_{234}^{2}
			\end{bsmallmatrix}.
\end{align*}
The obstructions from \cref{eq:tensor_FST,eq:tensor_SND} read
\begin{enumerate}[label=(\alph*)]
\item \label{it:osc_hyp_zer_1} $0 = -\frac{1}{2} \, H_{123} H_{134}$, %
\item \label{it:osc_hyp_zer_2} $0 = -\frac{1}{2} \, H_{123} H_{234}$, %
\item \label{it:osc_hyp_zer_3} $0 = \frac{1}{2} \, H_{123} H_{124}$, %
\item \label{it:osc_hyp_zer_4} $0 = -\frac{1}{2} \, H_{134} H_{234} a$, %
\item \label{it:osc_hyp_zer_5} $0 = \frac{1}{2} \, H_{124} H_{234} b$, %
\item \label{it:osc_hyp_zer_6} $0 = \frac{1}{2} \, H_{124} H_{134} (a - b)$, %
\item $0 = H_{134}^{2} a + a b {\lambda} + \frac{1}{2} \, H_{123}^{2} + \frac{3}{2} \, a^{2} + a b - \frac{1}{2} \, b^{2}$, %
\item $0 = H_{124}^{2} b + a b {\lambda} + \frac{1}{2} \, H_{123}^{2} - \frac{1}{2} \, a^{2} + a b + \frac{3}{2} \, b^{2}$, %
\item $0 = a b {\lambda} + H_{123}^{2} + \frac{1}{2} \, H_{234}^{2} + \frac{3}{2}$, %
\end{enumerate}
for the $\operatorname{FST}$-tensor and
\begin{enumerate}[label=(\alph*), resume]
\item $0 = -\frac{3}{2} \, H_{124} H_{134}^{2} a - \frac{3}{2} \, H_{124}^{3} b - H_{124} a b {\lambda} - \frac{3}{2} \, \cancel{H_{123}^{2} H_{124}} - \frac{3}{2} \, \cancel{H_{124} H_{234}^{2}} - H_{124} a b - 2 \, H_{124} b^{2}$, %
\item $0 = -\frac{3}{2} \, \cancel{H_{123} H_{134}^{2} a} - \frac{3}{2} \, \cancel{H_{123} H_{124}^{2} b} - H_{123} a b {\lambda} - \frac{3}{2} \, H_{123}^{3} - \frac{3}{2} \, \cancel{H_{123} H_{234}^{2}} - \frac{1}{2} \, H_{123} a^{2} - H_{123} a b - \frac{1}{2} \, H_{123} b^{2} - 2 \, H_{123}$, %
\item $0 = -\frac{3}{2} \, H_{134}^{3} a - \frac{3}{2} \, H_{124}^{2} H_{134} b - H_{134} a b {\lambda} - \frac{3}{2} \, \cancel{H_{123}^{2} H_{134}} - \frac{3}{2} \, \cancel{H_{134} H_{234}^{2}} - 2 \, H_{134} a^{2} - H_{134} a b$, %
\item $0 = -\frac{3}{2} \, \cancel{H_{134}^{2} H_{234} a} - \frac{3}{2} \, \cancel{H_{124}^{2} H_{234} b} - H_{234} a b {\lambda} - \frac{3}{2} \, \cancel{H_{123}^{2} H_{234}} - \frac{3}{2} \, H_{234}^{3} - \frac{1}{2} \, H_{234}$, %
\end{enumerate}
for the $\operatorname{SND}$-tensor.

We will now go through all possible cases and show that the only solutions are those claimed in the theorem.

\subparagraph{First case}

First assume that $H_{234} \ne 0$.
Then the conditions \labelcref{it:osc_hyp_zer_1,it:osc_hyp_zer_2,it:osc_hyp_zer_3,it:osc_hyp_zer_4,it:osc_hyp_zer_5,it:osc_hyp_zer_6} directly imply $H_{123} = H_{124} = H_{134} = 0$.
The equations then become
\begin{enumerate}
\item $0 = a b {\lambda} - \frac{1}{2} \, a^{2} + a b + \frac{3}{2} \, b^{2}$, %
\item $0 = -H_{234} a b {\lambda} - \frac{3}{2} \, H_{234}^{3} - \frac{1}{2} \, H_{234}$, %
\item $0 = a b {\lambda} + \frac{1}{2} \, H_{234}^{2} + \frac{3}{2}$ and %
\item $0 = a b {\lambda} + \frac{3}{2} \, a^{2} + a b - \frac{1}{2} \, b^{2}$. %
\end{enumerate}
The first minus the last equation then yield $0 = -2 a^2 + 2 b^2$, so $a = b$.
The last equation then shows $\lambda = -2$.
Finally, $\frac{1}{H_{234}}\, \mathrm{(ii)} \pm \mathrm{(iii)}$ yields $a = 1$ and $H_{234} = \pm 1$.
In other words,
\begin{equation*}
	[g_{ij}] = [\delta_{ij}]
	\quad\text{and}\quad
	H = \pm p^* \wedge q^* \wedge z^*.
\end{equation*}

\subparagraph{Second case}

Now assume $H_{234} = 0$ and $H_{123} \ne 0$.
Then $H_{124} = H_{134} = 0$ by \labelcref{it:osc_hyp_zer_1,it:osc_hyp_zer_3}.
The equations then simplify to
\begin{enumerate}
\item $0 = a b {\lambda} + H_{123}^{2} + \frac{3}{2}$, %
\item $0 = a b {\lambda} + \frac{1}{2} \, H_{123}^{2} - \frac{1}{2} \, a^{2} + a b + \frac{3}{2} \, b^{2}$, %
\item $0 = a b {\lambda} + \frac{1}{2} \, H_{123}^{2} + \frac{3}{2} \, a^{2} + a b - \frac{1}{2} \, b^{2}$ and %
\item $0 = -H_{123} a b {\lambda} - \frac{3}{2} \, H_{123}^{3} - \frac{1}{2} \, H_{123} a^{2} - H_{123} a b - \frac{1}{2} \, H_{123} b^{2} - 2 \, H_{123}$. %
\end{enumerate}

Again, the second minus the third equation shows $a = b$.
From here, the second equation shows $0 = a^{2} {\lambda} + \frac{1}{2} \, H_{123}^{2} + 2 \, a^{2}$, which leads to a contradiction, as
the last equation becomes
\begin{align*}
	0
	&= -H_{123} a b {\lambda} - \frac{3}{2} \, H_{123}^{3} - \frac{1}{2} \, H_{123} a^{2} - H_{123} a b - \frac{1}{2} \, H_{123} b^{2} - 2 \, H_{123} \\
	&= - H_{123} \pa[\bigg]{\underbrace{a^2 \lambda + \frac{1}{2} H_{123}^2 + 2a^2}_{= 0}\mathclose{} + H_{123}^2 + 2}
\end{align*}
and, by assumption, $H_{123}$ must not vanish.

\subparagraph{Third case}

Finally, assume $H_{234} = 0$ and $H_{123} = 0$.
Then
\begin{enumerate}
\item $0 = a b {\lambda} + \frac{3}{2}$, %
\item $0 = \frac{1}{2} \, H_{124} H_{134} (a - b)$, %
\item $0 = H_{124}^{2} b + a b {\lambda} - \frac{1}{2} \, a^{2} + a b + \frac{3}{2} \, b^{2}$, %
\item $0 = H_{134}^{2} a + a b {\lambda} + \frac{3}{2} \, a^{2} + a b - \frac{1}{2} \, b^{2}$, %
\item $0 = -H_{124} \pa*{\frac{3}{2} \, H_{134}^{2} a + \frac{3}{2} \, H_{124}^{2} b + a b {\lambda} + a b + 2 \, b^{2}}$ and %
\item $0 = -H_{134} \pa*{\frac{3}{2} \, H_{134}^{2} a + \frac{3}{2} \, H_{124}^{2} b + a b {\lambda} + 2 \, a^{2} + a b}$. %
\end{enumerate}

By substituting equation~(iii) into equation~(v), we obtain
\begin{align*}
	0
	&= -H_{124} \pa*{\frac{3}{2} \, H_{134}^{2} a + \frac{3}{2} \, H_{124}^{2} b + a b {\lambda} + a b + 2 \, b^{2}} \\
	&= -H_{124} \pa[\bigg]{\frac{3}{2} \, H_{134}^{2} a + \frac{1}{2} \, H_{124}^{2} b + \underbrace{H_{124}^{2} b + a b {\lambda} + a b + \frac{3}{2} \, b^{2}}_{= \frac{1}{2}\, a^2}\mathclose{} + \frac{1}{2} \, b^{2}} \\
	&= -H_{124} \underbrace{\pa[\bigg]{\frac{3}{2} \, H_{134}^{2} a + \frac{1}{2} \, H_{124}^{2} b + \frac{1}{2}\, a^2 + \frac{1}{2} \, b^{2}}}_{> 0}.
\end{align*}
Similarly, substituting equation~(iv) into equation~(vi), we get
\begin{align*}
	0
	&= -H_{134} \pa*{\frac{3}{2} \, H_{134}^{2} a + \frac{3}{2} \, H_{124}^{2} b + a b {\lambda} + 2 \, a^{2} + a b} \\
	&= -H_{134} \pa[\bigg]{\frac{1}{2} \, H_{134}^{2} a + \frac{3}{2} \, H_{124}^{2} b + \underbrace{H_{134}^{2} a + a b {\lambda} + \frac{3}{2} \, a^{2} + a b}_{= \frac{1}{2}\, b^2}\mathclose{} + \frac{1}{2} \, a^{2}} \\
	&= -H_{134} \underbrace{\pa[\bigg]{\frac{1}{2} \, H_{134}^{2} a + \frac{3}{2} \, H_{124}^{2} b + \frac{1}{2}\, b^2 + \frac{1}{2} \, a^{2}}}_{> 0}.
\end{align*}
Hence, $H = 0$.
The remaining, non-vanishing obstructions, (i), (iii) and (iv) are now readily solved and yield $\lambda = -2$ and $a = b = \frac{1}{2} \sqrt{3}$.
\end{proof}

\appendix

\section{Appendix}

\begin{proposition} \label{prop:diagonalization_hyperbolic_2x2}
Let $A$ denote a real, positive definite, symmetric $2 \times 2$-matrix.
Then there exists a matrix~$S \in \SO(1, 1)$, such that $S^\intercal A S$ is diagonal.
\end{proposition}
\begin{proof}
If $A$ is not already diagonal, without loss of generality, we can rescale $A$, so that
\begin{equation*}
    A
    =
    	\begin{bmatrix}
    		a & \pm 1 \\
    		\pm 1 & b
    	\end{bmatrix}.
\end{equation*}
By the principal minor criterion, $a > 0$, $b > 0$ and $ab - 1 > 0$, so the AM-GM-inequality yields
\begin{equation} \label{eq:am_gm_ineq}
    \frac{1}{2}(a + b) \ge \sqrt{ab} > 1.
\end{equation}
Now, the conjugation by a generic matrix in $\SO(1, 1)$ reads
\begin{equation*}
	\begin{bmatrix}
		x & y \\
		y & x
	\end{bmatrix}^\intercal
	\begin{bmatrix}
   		a & \pm 1 \\
   		\pm 1 & b
   	\end{bmatrix}
 	\begin{bmatrix}
 		x & y \\
 		y & x
 	\end{bmatrix}
 	=
 		\begin{bmatrix}
 			* & * \\
 			(a + b) xy \pm (x^2 + y^2) & *
 		\end{bmatrix},
\end{equation*}
where $x^2 - y^2 = 1$.
Thus, in order for this expression to be diagonal,
\begin{equation*}
	(a + b) xy \pm (x^2 + y^2) = 0.
\end{equation*}
By substituting $x = \pm \sqrt{1 + y^2}$ into this equation, we may get rid of the minus sign and it hence suffices to show the existence of some $y_0 \in \R$, such that the function
\begin{equation*}
	f(y) \ceq (a + b) \sqrt{1 + y^2} \cdot y + (1 + 2y^2)
\end{equation*}
vanishes in $y_0$.
Clearly, $f(0) = 1$ and
\begin{equation*}
	f(-\abs{y}) \le -(a + b) y^2 + 1 + 2y^2 = (2 - (a + b)) y^2 + 1,
\end{equation*}
so by \cref{eq:am_gm_ineq}, $f(y) \to -\infty$ for $y \to -\infty$.
In particular, we can find some value for which $f(y) < 0$, so by the intermediate value theorem, there exists a $y_0$ with $f(y_0) = 0$.
\end{proof}

\printbibliography

\end{document}